# LIKELIHOOD RATIO TESTS AND SINGULARITIES[1]

By Mathias Drton

*University of Chicago*

Many statistical hypotheses can be formulated in terms of polynomial equalities and inequalities in the unknown parameters and thus correspond to semi-algebraic subsets of the parameter space. We consider large sample asymptotics for the likelihood ratio test of such hypotheses in models that satisfy standard probabilistic regularity conditions. We show that the assumptions of Chernoff's theorem hold for semi-algebraic sets such that the asymptotics are determined by the tangent cone at the true parameter point. At boundary points or singularities, the tangent cone need not be a linear space and limiting distributions other than chi-square distributions may arise. While boundary points often lead to mixtures of chi-square distributions, singularities give rise to nonstandard limits. We demonstrate that minima of chi-square random variables are important for locally identifiable models, and in a study of the factor analysis model with one factor, we reveal connections to eigenvalues of Wishart matrices.

**1. Introduction.** Let $\mathcal{P}_\Theta = (P_\theta \mid \theta \in \Theta)$ be a parametric family of probability distributions on some measurable space. Suppose that $\Theta$ is an open subset of $\mathbb{R}^k$. For a hypothesis $\Theta_0 \subseteq \Theta$, consider testing

$$(1.1) \qquad H_0 : \theta \in \Theta_0 \quad \text{vs.} \quad H_1 : \theta \in \Theta \setminus \Theta_0$$

based on a large sample taken from a distribution in $\mathcal{P}_\Theta$. Under regularity conditions, the null distribution of the likelihood ratio statistic for the testing problem (1.1) can be approximated by the chi-square distribution $\chi_c^2$ with degrees of freedom $c$ equal to the codimension of $\Theta_0$, that is, $c = k - \dim(\Theta_0)$. The necessary regularity conditions combine probabilistic conditions on $\mathcal{P}_\Theta$ with geometric smoothness assumptions about $\Theta_0$. For example, the asymptotic approximation for the likelihood ratio test is valid when $\mathcal{P}_\Theta$ is a regular

Received March 2007; revised September 2007.
[1]Supported in part by NSF Grant DMS-05-05612 and the Institute for Mathematics and its Applications.
*AMS 2000 subject classifications.* 60E05, 62H10.
*Key words and phrases.* Algebraic statistics, factor analysis, large sample asymptotics, semi-algebraic set, tangent cone.

exponential family and $\Theta_0$ a smooth manifold, in which case the submodel $\mathcal{P}_{\Theta_0} = (P_\theta \mid \theta \in \Theta_0)$ is called a *curved exponential family* [18].

In this paper we consider the situation where probabilistic regularity conditions about $\mathcal{P}_\Theta$ hold but the null hypothesis $\Theta_0$ is a semi-algebraic set, that is, a set defined by polynomial equalities and inequalities in $\theta$. A semi-algebraic set has nice local geometric properties but it may have boundary points as well as singularities at which $\chi^2$-asymptotics are no longer valid. (For a rigorous definition of singularities of algebraic sets see, e.g., [3], Section 3.2 or [7], Section 9.) The case of semi-algebraic sets is important because many statistical hypotheses exhibit this special structure [9, 13]. Moreover, tools from algebraic geometry help in studying semi-algebraic sets and allow to make progress in the understanding of the likelihood ratio test.

Boundary points of statistical hypotheses have been discussed in the literature and often lead to asymptotic distributions that are mixtures of $\chi^2$-distributions. Two classic examples where boundary issues arise are variance component models [20] and factor analysis [25]; see also [24]. Singularities, however, do not seem to have received as much attention. For example, the parameter spaces of factor analysis models, which we will take up later in this paper, contain singularities at which the asymptotic distribution of the likelihood ratio statistic is not a $\chi^2$-mixture.

Issues with singularities can be illustrated nicely for hypotheses about the mean vector of a bivariate normal distribution $\mathcal{N}_2(\mu, I)$ with the covariance matrix equal to the identity matrix $I$. For a closed set $\Theta_0 \subseteq \Theta := \mathbb{R}^2$, the likelihood ratio statistic $\lambda_n$ for testing (1.1) is equal to the product of the sample size $n$ and the squared Euclidean distance between the sample mean vector and $\Theta_0$. The following two examples demonstrate nonstandard asymptotics for $\lambda_n$; the connection to tangent cones is based on a result of Chernoff [6] that we will revisit in this paper.

EXAMPLE 1.1 (Nodal cubic). Let $\Theta_0 = \{\mu \in \mathbb{R}^2 \mid \mu_2^2 = \mu_1^3 + \mu_1^2\}$ be the curve on the left in Figure 1, which can be parametrized as $\mathbf{f}(t) = [t^2 - 1, t(t^2-1)]$. The curve has a singularity at the point of self-intersection $\mu = 0$.

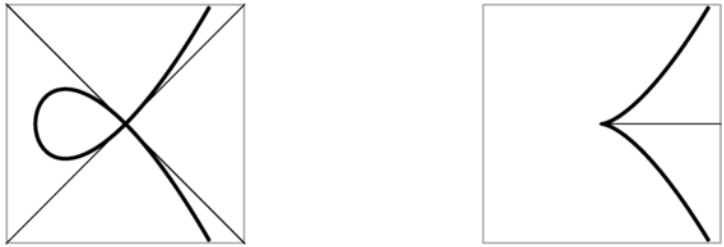

FIG. 1. *Nodal and cuspidal cubic.*

The lines $\mu_2 = \pm\mu_1$ in the plot indicate the tangent cone at $\mu = 0$. If $\mu = 0$ is the true parameter point and $n \to \infty$, then the likelihood ratio statistic $\lambda_n$ converges to the distribution of the squared Euclidean distance between a draw from $\mathcal{N}_2(0, I)$ and the lines $\mu_2 = \pm\mu_1$, that is, the distribution of the minimum of two independent $\chi_1^2$-random variables.

EXAMPLE 1.2 (Cuspidal cubic). Let $\Theta_0 = \{\mu \in \mathbb{R}^2 \mid \mu_2^2 = \mu_1^3\}$ be the curve with parametrization $\mathbf{f}(t) = (t^2, t^3)$ shown on the right in Figure 1. If the true parameter point is the cusp $\mu = 0$, then the asymptotic distribution of $\lambda_n$ is the mixture $\frac{1}{2}\chi_1^2 + \frac{1}{2}\chi_2^2$. This is the distribution of the squared Euclidean distance between a draw from $\mathcal{N}_2(0, I)$ and the tangent cone $\{\mu \mid \mu_1 \geq 0, \mu_2 = 0\}$.

In the above examples, points other than the origin are smooth points at which the curves locally look like a line. Thus, away from the origin the standard $\chi^2$-asymptotics, here $\chi_1^2$, apply. However, while $\chi^2$-limits arise almost everywhere, the convergence is not uniform and a very large sample size may be required for the $\chi^2$-distribution to provide a good approximation to the distribution of $\lambda_n$ if the true parameter is close to the singular locus. An important point is also that limiting distributions at singularities can be stochastically larger (Example 1.2) as well as smaller (Example 1.1) than the $\chi^2$-distribution obtained at smooth points.

The remainder of this paper begins with a review of the asymptotic theory for the likelihood ratio test (Section 2). We then show that the geometric regularity conditions in this theory are satisfied for semi-algebraic hypotheses (Section 3). In Section 4 we discuss algebraic methods that are helpful for determining tangent cones of semi-algebraic sets and can be used in particular to bound the asymptotic $p$-value of the likelihood ratio test. These methods are applied to factor analysis in Sections 5 and 6. Concluding remarks are given in Section 7.

**2. Likelihood ratio tests and tangent cones.** Suppose we observe a sample of independent and identically distributed random vectors $X^{(1)}, \ldots, X^{(n)} \in \mathbb{R}^m$ and that the distribution of $X^{(i)}$ belongs to the statistical model $\mathcal{P}_\Theta = (P_\theta \mid \theta \in \Theta)$. We assume that the distributions in $\mathcal{P}_\Theta$ are dominated by a common $\sigma$-finite measure $\nu$ with respect to which they have probability density functions $p_\theta : \mathbb{R}^m \to [0, \infty)$. For sample realizations $x^{(1)}, \ldots, x^{(n)}$, let

$$\ell_n : \Theta \to \mathbb{R}, \qquad \theta \mapsto \sum_{i=1}^n \log p_\theta(x^{(i)})$$

be the log-likelihood function of the model $\mathcal{P}_\Theta$. For $\Theta_0 \subseteq \Theta$, a *maximum likelihood estimator* $\hat{\theta}_{n,\Theta_0}$ in the (sub-)model $\mathcal{P}_{\Theta_0} = (P_\theta \mid \theta \in \Theta_0)$ satisfies

$$\ell_n(\hat{\theta}_{n,\Theta_0}) = \max_{\theta \in \Theta_0} \ell_n(\theta).$$

The *likelihood ratio statistic* for testing the fit of $\mathcal{P}_{\Theta_0}$, that is, for testing (1.1), is

$$(2.1) \qquad \lambda_n = 2\left(\sup_{\theta \in \Theta} \ell_n(\theta) - \sup_{\theta \in \Theta_0} \ell_n(\theta)\right).$$

For our study of large sample asymptotics for $\lambda_n$, we base ourselves on van der Vaart [34], Chapter 16, and make the following probabilistic regularity assumptions. Recall that the model $\mathcal{P}_\Theta$ is *differentiable in quadratic mean* at $\theta \in \Theta$ if there exists a measurable map $\dot{\ell}_\theta : \mathbb{R}^m \to \mathbb{R}^k$ such that

$$\lim_{h \to 0} \frac{1}{\|h\|^2} \int_{\mathbb{R}^m} \left(\sqrt{p_{\theta+h}(x)} - \sqrt{p_\theta(x)} - \frac{1}{2} h^t \dot{\ell}_\theta(x) \sqrt{p_\theta(x)}\right)^2 d\nu(x) = 0.$$

Lemma 7.6 in [34] gives a simple sufficient condition for differentiability in quadratic mean.

DEFINITION 2.1. A statistical model $\mathcal{P}_\Theta$ is *regular* at $\theta \in \Theta \subseteq \mathbb{R}^k$ if the following conditions hold:

(i) the point $\theta$ is in the interior of $\Theta$, which is assumed to be nonempty;

(ii) the model $\mathcal{P}_\Theta$ is differentiable in quadratic mean at $\theta$ with an invertible Fisher-information matrix $I(\theta) = \mathbb{E}_\theta[\dot{\ell}_\theta(X) \dot{\ell}_\theta(X)^t]$;

(iii) there exists a neighborhood $U(\theta) \subseteq \Theta$ of $\theta$ and a measurable function $\dot{\ell} : \mathbb{R}^m \to \mathbb{R}$, square-integrable as $\int_{\mathbb{R}^m} \dot{\ell}(x)^2 \, dP_\theta(x) < \infty$, such that

$$|\log p_{\theta_1}(x) - \log p_{\theta_2}(x)| \leq \dot{\ell}(x) \|\theta_1 - \theta_2\| \qquad \forall \theta_1, \theta_2 \in U(\theta).$$

(iv) the maximum likelihood estimator $\hat{\theta}_{n,\Theta}$ is consistent under $P_\theta$.

EXAMPLE 2.2. Let $\Theta = \mathbb{R}^{m \times m}_{\mathrm{pd}}$ be the cone of symmetric positive definite $m \times m$-matrices. The centered multivariate normal distributions $(\mathcal{N}_m(0, \Sigma) \mid \Sigma \in \Theta)$ form a regular exponential family (the natural parameter space is an open set). Such a family is regular in the sense of Definition 2.1 at all of its parameter points. In subsequent examples we tacitly identify the space of symmetric $m \times m$-matrices, denoted $\mathbb{R}^{m \times m}_{\mathrm{sym}}$, with $\mathbb{R}^{\binom{m+1}{2}}$ and index the vectors in the latter space by ordered pairs $ij$ with $i \leq j$. The inverse $I(\Sigma)^{-1}$ of the Fisher-information for $\Sigma = (\sigma_{ij})$ is then the $\binom{m+1}{2} \times \binom{m+1}{2}$-matrix with $(ij, k\ell)$-entry $\sigma_{ik}\sigma_{j\ell} + \sigma_{i\ell}\sigma_{jk}$.

For well-behaved large sample asymptotics of the likelihood ratio statistic $\lambda_n$ at a true parameter point $\theta_0$ in the null hypothesis $\Theta_0$, the probabilistic assumptions made above need to be complemented with assumptions about the local geometry of $\Theta_0$ at $\theta_0$. This local geometry expresses itself in the tangent cone.

DEFINITION 2.3. The tangent cone $T_\Theta(\theta)$ of the set $\Theta \subseteq \mathbb{R}^k$ at the point $\theta \in \mathbb{R}^k$ is the set of vectors in $\mathbb{R}^k$ that are limits of sequences $\alpha_n(\theta_n - \theta)$, where $\alpha_n$ are positive reals and $\theta_n \in \Theta$ converge to $\theta$.

The tangent cone $T_\Theta(\theta)$ is a closed set that is a *cone* in the sense that if $\tau \in T_\Theta(\theta)$ then the half-ray $\{\lambda\tau \mid \lambda \geq 0\}$ is contained in $T_\Theta(\theta)$. We state some properties of tangent cones that can be found for example in [23].

LEMMA 2.4. *Let $\theta \in \mathbb{R}^k$ and $\Theta, \Theta_1, \ldots, \Theta_m \subseteq \mathbb{R}^k$.*

(i) $T_{\Theta_1 \cup \cdots \cup \Theta_m}(\theta) = T_{\Theta_1}(\theta) \cup \cdots \cup T_{\Theta_m}(\theta)$.
(ii) $T_{\Theta_1 \cap \cdots \cap \Theta_m}(\theta) \subseteq T_{\Theta_1}(\theta) \cap \cdots \cap T_{\Theta_m}(\theta)$.
(iii) *If $\Theta - \theta$ is a cone, then $T_\Theta(\theta)$ is the topological closure of $\Theta - \theta$.*
(iv) *Let $\Theta = \mathbf{f}(\Gamma)$ for some differentiable map $\mathbf{f}: \mathbb{R}^d \to \mathbb{R}^k$ and some open set $\Gamma \subseteq \mathbb{R}^d$. If $\theta = \mathbf{f}(\gamma)$ for some $\gamma \in \Gamma$, then $T_\Theta(\theta)$ contains the linear space spanned by the columns of the Jacobian*

$$\left(\frac{\partial \mathbf{f}_i(\gamma)}{\partial \gamma_j}\right) \in \mathbb{R}^{k \times d}.$$

The next definition describes a regularity requirement on the local geometry of a hypothesis $\Theta_0$ at a point $\theta_0$; see [23], Proposition 6.2 and [14].

DEFINITION 2.5. The set $\Theta \subseteq \mathbb{R}^k$ is Chernoff-regular at $\theta$ if for every vector $\tau$ in the tangent cone $T_\Theta(\theta)$ there exist $\varepsilon > 0$ and a map $\alpha:[0,\varepsilon) \to \Theta$ with $\alpha(0) = \theta$ such that $\tau = \lim_{t \to 0+}[\alpha(t) - \alpha(0)]/t$.

Under Chernoff-regularity, likelihood ratio statistics converge to Mahalanobis distances between a random draw from a multivariate normal distribution and the tangent cone at the true parameter point $\theta_0$. This result first appeared in [6]. The version given here is proven in [34], Theorem 16.7. Note that under Chernoff-regularity the sets $\sqrt{n}(\Theta_0 - \theta_0)$ converge to $T_{\Theta_0}(\theta_0)$ in the sense of [34]; compare [14].

THEOREM 2.6. *Let $\theta_0 \in \Theta_0 \subseteq \Theta \subseteq \mathbb{R}^k$ be a true parameter point at which the model $\mathcal{P}_\Theta$ is regular with Fisher-information $I(\theta_0)$. Let $\bar{Z} \sim \mathcal{N}_k(0, I(\theta_0)^{-1})$. If $\Theta_0$ is Chernoff-regular at $\theta_0$ and the maximum likelihood estimator $\hat{\theta}_{n,\Theta_0}$ is consistent, then as $n$ tends to infinity, the likelihood ratio statistic $\lambda_n$ converges to the distribution of the squared Mahalanobis distance*

$$\min_{\tau \in T_{\Theta_0}(\theta_0)} (\bar{Z} - \tau)^t I(\theta_0)(\bar{Z} - \tau).$$

*If $I(\theta_0) = I(\theta_0)^{t/2} I(\theta_0)^{1/2}$ and $Z \sim \mathcal{N}_k(0, I)$ is a standard normal vector, then the squared Mahalanobis distance has the same distribution as the* *squared Euclidean distance*

$$\min_{\tau \in T_{\Theta_0}(\theta_0)} \|Z - I(\theta_0)^{1/2}\tau\|^2$$

*between $Z$ and the linearly transformed tangent cone $I(\theta_0)^{1/2}T_{\Theta_0}(\theta_0)$.*

We remark that changing the matrix square root $I(\theta_0)^{1/2}$ corresponds to an orthogonal transformation under which Euclidean distances as well as the standard normal distribution are invariant.

If the tangent cone $T_{\Theta_0}(\theta_0)$ in Theorem 2.6 is a $d$-dimensional linear subspace of $\mathbb{R}^k$, then we recover the standard $\chi^2$-asymptotics because the squared Euclidean distance between a $d$-dimensional subspace and a standard normal vector follows a $\chi^2_{k-d}$-distribution. Another well-known case arises if $T_{\Theta_0}(\theta_0)$ is a convex cone. In this case the limiting distribution is a mixture of $\chi^2$-distributions with degrees of freedom larger than or equal to the codimension of $T_{\Theta_0}(\theta_0)$; see [19, 25, 26, 28, 29, 31] and Example 1.2. These mixtures are also known as chi-bar-square distributions.

The next example gives another important type of nonstandard limiting distributions that we have already encountered in the artificial Example 1.1.

EXAMPLE 2.7. Suppose we wish to test whether the marginal independence $X_1 \perp\!\!\!\perp X_2$ and the conditional independence $X_1 \perp\!\!\!\perp X_2 \mid X_3$ hold simultaneously in $(X_1, X_2, X_3)^t \sim \mathcal{N}_3(0, \Sigma)$. This is the case if and only if the unknown covariance matrix $\Sigma = (\sigma_{ij})$ satisfies that $\sigma_{12} = 0$ and $\sigma_{13}\sigma_{23} = 0$. Define the linear space $L_i = \{z \in \mathbb{R}^{3\times 3}_{\text{sym}} \mid z_{12} = z_{i3} = 0\}$ for $i = 1, 2$. The null hypothesis $\Theta_0 \subseteq \mathbb{R}^{3\times 3}_{\text{sym}}$ comprises the positive definite matrices in $L_1 \cup L_2$. A true covariance matrix $\Sigma_0$ that is not diagonal belongs either to $L_1$ or to $L_2$ such that the tangent cone $T_{\Theta_0}(\Sigma_0)$ is equal to $L_1$ or $L_2$, respectively. Since $\dim(L_i) = 4$, it follows that $\lambda_n$ converges to a $\chi^2_{6-4} = \chi^2_2$-distribution. If, however, $\Sigma_0$ is diagonal then $\Theta_0 - \Sigma_0$ coincides with the closed cone $L_1 \cup L_2$ near the origin and, by Lemma 2.4(iii), $T_{\Theta_0}(\Sigma_0) = L_1 \cup L_2$. The Fisher-information $I(\Sigma_0)$ and its symmetric square-root $I(\Sigma_0)^{1/2}$ are now diagonal. Diagonal transformations leave the cone $L_1 \cup L_2$ invariant such that by Theorem 2.6, $\lambda_n$ converges in distribution to the minimum of $Z_{12}^2 + Z_{13}^2$ and $Z_{12}^2 + Z_{23}^2$ for a standard normal random vector $Z \in \mathbb{R}^{3\times 3}_{\text{sym}}$. This is the distribution of $W_{12} + \min(W_{13}, W_{23})$, where $W_{12}$, $W_{13}$ and $W_{23}$ are independent $\chi^2_1$-random variables. We note that this example can also be worked out by elementary means [8].

Examples in which Chernoff-regularity fails and the likelihood ratio statistic $\lambda_n$ does not converge in distribution can be found in [9] and [14].

REMARK 2.8 (Comparing nested submodels). The above setup considers problem (1.1) that compares the fit of the submodel $\mathcal{P}_{\Theta_0}$ to the fit of the saturated model $\mathcal{P}_\Theta$. Instead, one may wish to compare to the fit of another submodel $\mathcal{P}_{\Theta_1}$, that is, test $H_0: \theta \in \Theta_0$ versus $H_1: \theta \in \Theta_1 \setminus \Theta_0$, where $\Theta_0 \subseteq \Theta_1 \subseteq \Theta$. However, if Theorem 2.6 applies to both (1.1) and the problem $H_0: \theta \in \Theta_1$ versus $H_1: \theta \in \Theta \setminus \Theta_1$, then we can deduce that the asymptotic distribution of the likelihood ratio statistic

$$\lambda_n = 2\bigg(\sup_{\theta \in \Theta_1} \ell_n(\theta) - \sup_{\theta \in \Theta_0} \ell_n(\theta)\bigg)$$

is given by the difference of the squared Mahalanobis distances between the random vector $Z \sim \mathcal{N}_k(0, I(\theta_0)^{-1})$ and the two tangent cones, namely

$$\lambda_n \xrightarrow{n \to \infty}_d \bigg[\min_{\theta \in T_{\Theta_0}(\theta_0)} (Z-\theta)^t I(\theta_0)(Z-\theta)\bigg] - \bigg[\min_{\theta \in T_{\Theta_1}(\theta_0)} (Z-\theta)^t I(\theta_0)(Z-\theta)\bigg].$$

REMARK 2.9 (Maximum likelihood estimators). Under the conditions of Theorem 2.6, the maximum likelihood estimator $\hat{\theta}_{n,\Theta_0}$ can be shown to be distributed like the projection of $Z \sim \mathcal{N}_k(0, I(\theta_0)^{-1})$ on the tangent cone $T_{\Theta_0}(\theta_0)$ [34], Theorem 7.12. Here, projection refers to the minimizer of the Mahalanobis distance $(Z-\theta)^t I(\theta_0)(Z-\theta)$ over $\theta \in T_{\Theta_0}(\theta_0)$. This minimizer is almost surely unique [14], Proposition 4.2. For results on local maximizers of the likelihood function, see [27].

**3. Semi-algebraic hypotheses.** In principle, Chernoff regularity has to be verified in every application of Theorem 2.6. However, as we detail in this section, Chernoff regularity is automatic if the hypothesis $\Theta_0$ is given by a semi-algebraic set. The map $\alpha$ in Definition 2.5 can be chosen very smooth for semi-algebraic sets.

3.1. *Semi-algebraic sets.* We begin by briefly reviewing some of the properties of semi-algebraic sets. In-depth treatments can be found in [2, 3, 5].

DEFINITION 3.1. Let $\mathbb{R}[\mathbf{t}] = \mathbb{R}[t_1, \ldots, t_k]$ be the ring of polynomials in the indeterminates $t_1, \ldots, t_k$ with real coefficients. A *basic semi-algebraic set* is a subset of $\mathbb{R}^k$ that is of the form

$$\Theta = \{\theta \in \mathbb{R}^k \mid f(\theta) = 0 \; \forall f \in F, h(\theta) > 0 \; \forall h \in H\},$$

where $F, H \subset \mathbb{R}[\mathbf{t}]$ are (possibly empty) collections of polynomials and $H$ is finite. A *semi-algebraic set* is a finite union of basic semi-algebraic sets. If $H = \varnothing$ then $\Theta$ is called a *real algebraic variety*.

Complements, finite unions and finite intersections of semi-algebraic sets are again semi-algebraic. If $\Gamma \subseteq \mathbb{R}^d$ is semi-algebraic and $\mathbf{f} : \mathbb{R}^d \to \mathbb{R}^k$ a polynomial map (or a rational map defined everywhere on $\Gamma$), then the image $\mathbf{f}(\Gamma)$ is a semi-algebraic set. The parameter spaces of many statistical models are such images; compare [22].

A semi-algebraic set $\Theta$ can be written as a disjoint union of finitely many smooth manifolds $\Theta_1, \ldots, \Theta_s$ such that if $\Theta_i$ and the closure $\mathrm{cl}(\Theta_j)$ have a nonempty intersection then $\Theta_i \subseteq \mathrm{cl}(\Theta_j)$ and $\dim(\Theta_i) < \dim(\Theta_j)$. Such a partition is known as a stratification of $\Theta$. The *dimension* of $\Theta$ can be defined as the largest dimension of any manifold in the stratification. If $\Theta = \mathbf{f}(\Gamma)$ is the image of an open semi-algebraic set $\Gamma$ under a polynomial map $\mathbf{f}$, then $\dim(\Theta)$ is equal to the maximal rank of any Jacobian of $\mathbf{f}(\gamma)$ for $\gamma \in \Gamma$. This maximal rank is achieved at almost every $\gamma$.

DEFINITION 3.2. For a point $\theta$ in the semi-algebraic set $\Theta$, define $d_m$ to be the dimension of the semi-algebraic set $\Theta \cap B_{1/m}(\theta)$, where $B_{1/m}(\theta)$ is the open ball of radius $1/m$ around $\theta$. The sequence $(d_m)_{m \in \mathbb{N}}$ being nonincreasing, there exists $m_0$ such that $d_m = d_{m_0}$ for all $m \geq m_0$. The *local dimension* of $\theta$ is defined to be $\dim_\theta(\Theta) := d_{m_0}$, which is no larger than $\dim(\Theta)$.

If there exists a ball $B_r(\theta)$ such that the semi-algebraic set $\Theta \cap B_r(\theta)$ is a $d$-dimensional smooth manifold then $\theta$ is a *smooth point* and its local dimension is $\dim_\theta(\Theta) = d$. Almost all points $\theta$ of a semi-algebraic set $\Theta$ are smooth of local dimension $\dim_\theta(\Theta) = \dim(\Theta)$. In other words, the set of points $\theta \in \Theta$ with $\dim_\theta(\Theta) < \dim(\Theta)$ is a subset of dimension smaller than $\dim(\Theta)$.

A semi-algebraic set, even a real algebraic variety, may have smooth points of different local dimensions. For example, the so-called Whitney umbrella defined by $x^2 z = y^2$ in $\mathbb{R}^3$ has two-dimensional smooth points, which arise if $x \neq 0$ or $y \neq 0$. The points with $x = y = 0$ and $z \geq 0$ are not smooth, but the points $x = y = 0$ and $z < 0$ lie on a line and are thus one-dimensional smooth points. However, if $\Theta = \mathbf{f}(\Gamma)$ for an open semi-algebraic set $\Gamma$ and a polynomial map $\mathbf{f}$, then $\Theta$ is pure-dimensional in the sense that $\dim_\theta(\Theta) = \dim(\Theta)$ for all $\theta \in \Theta$ [11].

3.2. *Tangent cones of semi-algebraic sets.* If $\Theta$ is semi-algebraic, then the tangent cone $T_\Theta(\theta)$ at a point $\theta \in \Theta$ is also a semi-algebraic set. The dimension of the tangent cone $T_\Theta(\theta)$ is at most $\dim(\Theta)$ and may be strictly smaller even for polynomial images of open semi-algebraic sets. For example, if $\mathbf{f} : (s,t) \mapsto (s^2 + t^2, s^3 + t^3)$ then $\mathbf{f}(\mathbb{R}^2)$ is the two-dimensional set that includes all points that are on or to the right of the cuspidal cubic shown in Figure 1. At the origin, this two-dimensional set has the one-dimensional dashed half-ray as tangent cone.

Despite this possible difference between the dimension of the tangent cone and the dimension of the set itself, the tangent cones to semi-algebraic sets are very well-behaved: the vectors in the tangent cone are tangent to very smooth curves in the semi-algebraic set. This has the following important consequence that ensures the existence of limiting distributions in many examples.

LEMMA 3.3. *A semi-algebraic set $\Theta \subseteq \mathbb{R}^k$ is everywhere Chernoff regular.*

PROOF. Proposition 2 of [21] shows that if $\Theta$ is a real algebraic variety and $\tau \in T_\Theta(\theta)$ for some $\theta \in \Theta$, then there exists a real analytic curve $\alpha : [0, \varepsilon) \to \Theta$ with $\alpha(0) = \theta$ that for $t \in [0, \varepsilon)$ has a convergent Taylor series expansion as $\alpha(t) = \tau t^r + O(t^{r+1})$ with $r \geq 1$. The result for semi-algebraic sets can be proven in exactly the same fashion by altering Claim 2 in [21] to make a requirement about inequalities $h_j(\theta) > 0$. By Lemma 2.4(i) it suffices to consider a basic semi-algebraic set. □

REMARK 3.4. Let $(P_\theta \mid \theta \in \Theta)$ be a regular exponential family. Drton and Sullivant [9] define a submodel $(P_\theta \mid \theta \in \Theta_0)$ to be an *algebraic exponential family* if $\Theta_0 \subset \Theta$ is diffeomorphic to a semi-algebraic set. By Lemma 3.3, $\Theta_0$ is everywhere Chernoff-regular such that Theorem 2.6 applies at every point $\theta_0 \in \Theta_0$ at which the maximum likelihood estimator $\hat{\theta}_{n,\Theta_0}$ is consistent. According to [4], Theorem 3.1, Corollary 3.3, $\hat{\theta}_{n,\Theta_0}$ is consistent if $\Theta_0$ is locally compact at $\theta_0$. A semi-algebraic set need not be locally compact. However, the likelihood ratio statistic $\lambda_n$ does not change if in (2.1) we replace $\Theta_0$ by the union of $\Theta_0$ and the closure of $B_\epsilon(\theta_0) \cap \Theta_0$. This closure is meaningful for small $\epsilon$ because $\theta_0$ is in the interior of $\Theta$. With this change $\Theta_0$ is locally compact at $\theta_0$, and we can deduce that the first-order asymptotics of the likelihood ratio test for testing the goodness-of-fit of an algebraic exponential family are always given by Mahalanobis distances from the tangent cone.

3.3. *Locally identifiable models.* Suppose $\mathcal{P}_\Theta = (P_\theta \mid \theta \in \Theta)$ is an identifiable model, that is, $P_\theta = P_{\bar{\theta}}$ implies $\theta = \bar{\theta}$. Let $\Gamma \subseteq \mathbb{R}^d$ be an open semi-algebraic set and $\mathbf{f} : \Gamma \to \mathbb{R}^k$ a polynomial or rational map such that $\Theta_0 = \mathbf{f}(\Gamma) \subseteq \Theta$. The submodel $\mathcal{P}_{\Theta_0}$ that is parameterized by $\mathbf{f}$ is globally identifiable at $\gamma_0 \in \Gamma$ if $\gamma_0$ is the unique point in $\Gamma$ that is mapped to $\theta_0 = \mathbf{f}(\gamma_0)$. The submodel $\mathcal{P}_{\Theta_0}$ is *locally identifiable* at $\gamma_0 \in \Gamma$ if there exists a neighborhood $U(\gamma_0) \subseteq \Gamma$ of $\gamma_0$ such that $\mathbf{f}^{-1}(\theta_0) \cap U(\gamma_0) = \{\gamma_0\}$. Local identifiability often arises as assumed in the following proposition, where $J(\gamma)$ denotes the Jacobian of $\mathbf{f}$ at $\gamma$.

PROPOSITION 3.5. *Let $\theta_0$ be a true parameter point at which $\mathcal{P}_\Theta$ is regular and the maximum likelihood estimator $\hat{\theta}_{n,\Theta_0}$ consistent. Suppose $\mathbf{f}^{-1}(\theta_0)$ is a finite set such that the Jacobian $J(\gamma)$ has full rank at all $\gamma \in \mathbf{f}^{-1}(\theta_0)$. If $\mathbf{f}$ is proper at $\theta_0$, that is, there exists a compact neighborhood $V \subseteq \mathbb{R}^k$ of $\theta_0$ such that $\mathbf{f}^{-1}(V) \cap \Gamma$ is compact in $\mathbb{R}^d$, then the likelihood ratio statistic $\lambda_n$ for (1.1) converges to the distribution of a minimum of at most $|\mathbf{f}^{-1}(\theta_0)|$ random variables with $\chi^2_{k-\dim(\Theta_0)}$-distribution.*

PROOF. By the full rank assumption, there exist neighborhoods $U(\gamma) \subseteq \Gamma$ of $\gamma \in \mathbf{f}^{-1}(\theta_0)$ such that $M_\gamma = \mathbf{f}(U(\gamma))$ are smooth manifolds. Consider a sequence $(\theta_n) = (\mathbf{f}(\gamma_n))$ in $\Theta_0$ that converges to $\theta_0$. Since $\theta_n \in V$ for all large $n$, the sequence $(\gamma_n)$ is eventually contained in the compactum $\mathbf{f}^{-1}(V) \cap \Gamma$. Since $\mathbf{f}$ is continuous, all accumulation points of $(\gamma_n)$ are in the finite preimage $\mathbf{f}^{-1}(\theta_0)$. Therefore,

$$\gamma_n \in \bigcup_{\gamma \in \mathbf{f}^{-1}(\theta_0)} U(\gamma)$$

for all $n$ larger than some $n_0 \in \mathbb{N}$. It follows that locally at $\theta_0$, the set $\Theta_0$ is equal to a finite union of the smooth manifolds $M_\gamma$, $\gamma \in \mathbf{f}^{-1}(\theta_0)$. According to Lemma 2.4(i), the tangent cone $T_{\Theta_0}(\theta_0)$ is the finite union of the tangent spaces of the manifolds $M_\gamma$, which are the linear spaces $L_\gamma$ spanned by the columns of $J(\gamma)$.

Let $Z \sim \mathcal{N}_k(0, I(\theta_0))$. The limiting distribution of $\lambda_n$ is the distribution of

$$\min_{\gamma \in \mathbf{f}^{-1}(\theta_0)} \left( \min_{\tau \in L_\gamma} (Z - \tau)^t I(\theta_0)(Z - \tau) \right).$$

Since the $L_\gamma$ are linear spaces of dimension $\dim(\Theta_0)$, the displayed expression is a minimum of $\chi^2_{k-\dim(\Theta_0)}$-random variables. If $L_\gamma = L_{\bar\gamma}$ for $\gamma \neq \bar\gamma \in \mathbf{f}^{-1}(\theta_0)$, then only one of two associated $\chi^2$-variables needs to be included in the minimum. □

Proposition 3.5 makes no statement about the dependence of the $\chi^2$-random variables in the minimum. In the artificial Example 1.1, the two $\chi^2_1$-random variables were independent but, as illustrated next, this is not the case in general.

EXAMPLE 3.6. Suppose $\varepsilon_1, \ldots, \varepsilon_4$ are independent normal random variables distributed as $\varepsilon_i \sim \mathcal{N}(0, \omega_i)$ with $\omega_i > 0$. Consider the system of linear equations

$$Y_1 = \varepsilon_1, \qquad Y_2 = \beta_{21} Y_1 + \beta_{24} Y_4 + \varepsilon_2,$$
$$Y_3 = \beta_{31} Y_1 + \beta_{32} Y_2 + \varepsilon_3, \qquad Y_4 = \beta_{43} Y_3 + \varepsilon_4,$$

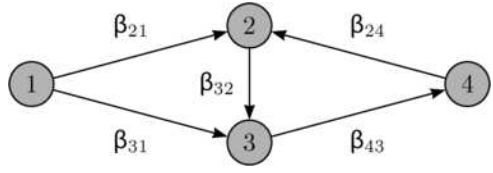

FIG. 2. *Graph with feedback loop.*

which has a feedback loop among $Y_2$, $Y_3$ and $Y_4$; compare the graphical representation in Figure 2. In matrix form, the equations state that if $Y = (Y_1, \ldots, Y_4)^t$ and $\varepsilon = (\varepsilon_1, \ldots, \varepsilon_4)^t$ then $BY = \varepsilon$ for

$$(3.1) \qquad B = \begin{pmatrix} 1 & 0 & 0 & 0 \\ -\beta_{21} & 1 & 0 & -\beta_{24} \\ -\beta_{31} & -\beta_{32} & 1 & 0 \\ 0 & 0 & -\beta_{43} & 1 \end{pmatrix}.$$

Let

$$D = \{(\beta_{21}, \ldots, \beta_{43})^t \in \mathbb{R}^5 \mid \det(B) = 1 - \beta_{32}\beta_{24}\beta_{43} \neq 0\}$$

and $\Gamma = D \times (0, \infty)^4$. The map from $(\beta, \omega) \in \Gamma$ to the covariance matrix of $Y$ is rational; denote it by $\mathbf{f}$. The set $\Theta_0 = \mathbf{f}(\Gamma)$ is a semi-algebraic subset of the cone of positive definite matrices $\Theta = \mathbb{R}^{4 \times 4}_{\text{pd}}$. It is the parameter space of the Gaussian model $\mathcal{P}_{\Theta_0} = (\mathcal{N}_4(0, \Sigma) \mid \Sigma \in \Theta_0)$ that is induced by the equation system.

The Jacobian of $\mathbf{f}$ is of rank at least 8 for all $(\beta, \omega) \in \Gamma$. It is of full rank 9 unless

$$(3.2) \qquad \beta_{31} + \beta_{32}\beta_{21} = 0 \quad \text{and} \quad \beta_{32}\beta_{43}\beta_{24} = -1.$$

Details of calculations that yield this and other facts employed here are given in Appendix A.1. By Lemma A.1, the model $\mathcal{P}_{\Theta_0}$ is globally identifiable at $(\beta, \omega)$ unless

$$(3.3) \quad \beta_{31} + \beta_{32}\beta_{21} = 0, \qquad \beta_{32}\beta_{43}\beta_{24} \neq -1 \quad \text{and} \quad \beta_{32}, \beta_{43}, \beta_{24} \neq 0.$$

If $\beta$ satisfies (3.3), then $\mathcal{P}_{\Theta_0}$ is locally identifiable with the preimage of $\Sigma = \mathbf{f}(\beta, \omega)$ always being of cardinality two. Moreover, by Lemma A.2, the map $\mathbf{f}$ is proper at $\Sigma$. Hence, in this locally identifiable case, the likelihood ratio statistic $\lambda_n$ converges to the distribution of the minimum of two $\chi^2_1$-random variables.

Suppose the true parameter point is $\Sigma_0 = \mathbf{f}(\beta, \omega)$ with $\beta$ as in (3.3). Using Lemma 2.4(iv), it can be shown that the tangent cone $T_{\Theta_0}(\Sigma_0)$ is equal to the union of two hyperplanes whose normal vectors $\eta$ and $\bar{\eta}$ have zero components except at their 13- and 14-entries. The nontrivial entries are

$$(3.4) \qquad \eta_{(13,14)} = (\beta_{43}, -1)^t$$

and

$$\bar{\eta}_{(13,14)} = \left(1, -\frac{\beta_{43}(\omega_3 + \beta_{32}^2\omega_2 + \beta_{32}^2\beta_{24}^2\omega_4)}{\omega_4 + \beta_{43}^2\omega_3 + \beta_{32}^2\beta_{43}^2\omega_2}\right)^t. \tag{3.5}$$

Equation (3.4) is readily obtained by computing the kernel of the transposed Jacobian of $\mathbf{f}$ at $(\beta, \omega)$. Equation (3.5) follows from replacing $\beta_{43}$ by the component $\bar{\beta}_{43}$ of the second vector $(\bar{\beta}, \bar{\omega})$ with $\mathbf{f}(\bar{\beta}, \bar{\omega}) = \Sigma_0$; compare (A.2) where $\kappa_i = 1/\omega_i$.

In order to describe the limiting distribution of the likelihood ratio statistic more precisely, we need to consider the transformed tangent cone $I(\Sigma_0)^{1/2}T_{\Theta_0}(\theta_0)$, which is the union of the two hyperplanes with normal vectors $I(\Sigma_0)^{-t/2}\eta$ and $I(\Sigma_0)^{-t/2}\bar{\eta}$. The cosine of the angle between these two normal vectors is equal to

$$\rho = \frac{\eta_1^t I(\Sigma_0)^{-1}\eta_2}{\sqrt{\eta_1^t I(\Sigma_0)^{-1}\eta_1 \cdot \eta_2^t I(\Sigma_0)^{-1}\eta_2}}.$$

This expression simplifies to

$$\rho = \frac{\beta_{43}\omega_3 + \beta_{43}\beta_{32}^2\omega_2 - \beta_{24}\beta_{32}\omega_4}{\sqrt{(\omega_3 + \beta_{32}^2\omega_2 + \beta_{24}^2\beta_{32}^2\omega_4)(\omega_4 + \beta_{43}^2\omega_3 + \beta_{32}^2\beta_{43}^2\omega_2)}}; \tag{3.6}$$

recall the formula for the inverse Fisher-information matrix $I(\Sigma_0)^{-1}$ from Example 2.2. We may thus conclude that $\lambda_n$ converges to the distribution of the squared Euclidean distance between a standard normal point $Z \sim \mathcal{N}_2(0, I)$ in $\mathbb{R}^2$ and two lines through the origin that intersect at angle $\cos^{-1}(\rho)$.

If $\beta_{31} + \beta_{32}\beta_{21} = 0$ and at least one of the parameters $\beta_{32}$, $\beta_{43}$, $\beta_{24}$ zero, then $\mathcal{P}_{\Theta_0}$ is globally identifiable at $(\beta, \omega)$, the Jacobian of $\mathbf{f}$ of full rank, but Proposition 3.5 does not apply as explained at the end of Appendix A.1. It is interesting that in this case, the limiting distribution of $\lambda_n$ is not a $\chi^2$-distribution, which is shown in Proposition A.4 in Appendix A.2. This fact as well as results on the remaining cases are obtained using the algebraic techniques we present in the next section.

**4. Algebraic tangent cones and bounds on $p$-values.** In this section, we explain how algebraic tools can help find smooth and non-smooth points of a semi-algebraic hypothesis. Algebraic methods also allow one to compute (asymptotic) bounds on $p$-values for the likelihood ratio test.

4.1. *Boundary points and singularities.* For a semi-algebraic set $\Theta \subseteq \mathbb{R}^k$, define $\mathcal{I}(\Theta)$ to be the ideal of polynomials $f \in \mathbb{R}[\mathbf{t}]$ that vanish whenever evaluated at a point $\theta \in \Theta$. By Hilbert's basis theorem, the ideal $\mathcal{I}(\Theta)$ has a finite generating set $\{f_1, \ldots, f_\ell\} \subseteq \mathbb{R}[\mathbf{t}]$. In other words, there exist finitely many polynomials $f_1, \ldots, f_\ell$ such that $f \in \mathcal{I}(\Theta)$ if and only if $f = h_1 f_1 + \cdots + h_\ell f_\ell$ for some $h_j \in \mathbb{R}[\mathbf{t}]$. The real algebraic variety defined by the vanishing of all polynomials in $\mathcal{I}(\Theta)$, or equivalently the polynomials $f_1, \ldots, f_\ell$, is the *Zariski closure* $\bar{\Theta}$ of $\Theta$. The Zariski closure is the smallest real algebraic variety containing $\Theta$ and in particular $\dim(\bar{\Theta}) = \dim(\Theta)$.

DEFINITION 4.1. Let $\Theta \subseteq \mathbb{R}^k$ be a semi-algebraic set with Zariski closure $\bar{\Theta}$. A subset of $\Theta$ is open in $\bar{\Theta}$ if it is equal to the intersection of an open set in $\mathbb{R}^k$ and $\bar{\Theta}$. The interior $\mathrm{int}(\Theta)$ is the union of all subsets of $\Theta$ that are open in $\bar{\Theta}$. The boundary $\mathrm{bd}(\Theta)$ is the complement $\Theta \setminus \mathrm{int}(\Theta)$.

At a boundary point $\theta \in \mathrm{bd}(\Theta)$ the tangent cone need not be a linear space such that nonstandard limiting distributions may arise for the likelihood ratio statistic. However, this phenomenon may also occur at singularities. We recall the definition of singularities as given, for example, in [3], Section 3.2.

A real algebraic variety $\Theta$ is *irreducible* if it cannot be written as the union of two strict subsets that are also real algebraic varieties. Any real algebraic variety $\Theta$ can be written as a finite union of irreducible varieties,

$$\Theta = \Theta_1 \cup \cdots \cup \Theta_\ell. \tag{4.1}$$

If no two varieties $\Theta_i$ and $\Theta_j$ in (4.1) are ordered by inclusion, then the decomposition in (4.1) is unique up to reordering and the irreducible varieties $\Theta_i$ are called the *irreducible components* of $\Theta$. Let $\{f_1, \ldots, f_\ell\} \subseteq \mathbb{R}[\mathbf{t}]$ generate the ideal $\mathcal{I}(\Theta)$, and let $J(\theta) \in \mathbb{R}^{\ell \times k}$ be the Jacobian with $ij$th entry $\partial f_i(\theta)/\partial \theta_j$. Let $r(\Theta)$ be the maximum rank of any matrix $J(\theta), \theta \in \Theta$. Then $r(\Theta)$ is independent of the choice of the generating set $\{f_1, \ldots, f_\ell\}$ and it holds that $r(\Theta) = k - \dim(\Theta)$.

DEFINITION 4.2. Let $\theta$ be a point in the real algebraic variety $\Theta \subseteq \mathbb{R}^k$.

(i) If $\Theta$ is irreducible and the rank of $J(\theta)$ is smaller than $r(\Theta)$ then $\theta$ is a singular point of $\Theta$.

(ii) If $\Theta_1, \ldots, \Theta_\ell$ are the irreducible components of $\Theta$ then $\theta$ is a singular point of $\Theta$ if it is a singular point of some $\Theta_i$ or if it is in an intersection $\Theta_i \cap \Theta_j$.

If $\Theta$ is a semi-algebraic set then $\theta \in \Theta$ is a singular point of $\Theta$ if it is a singular point of the Zariski closure $\bar{\Theta}$.

The software `Singular` [15] provides routines for computing all singularities of $\Theta$ from a generating set $\{f_1, \ldots, f_\ell\} \subseteq \mathbb{R}[\mathbf{t}]$ for $\mathcal{I}(\Theta)$.

A nonsingular interior point $\theta$ of a semi-algebraic set $\Theta$ is also a smooth point with tangent cone $T_\Theta(\theta)$ that is equal to the linear space given by the kernel of the Jacobian matrix $J(\theta)$. This fact translates into the following statistical result.

THEOREM 4.3. *Let $\theta_0 \in \Theta_0 \subseteq \Theta \subseteq \mathbb{R}^k$ be a true parameter point at which the model $\mathcal{P}_\Theta$ is regular and the maximum likelihood estimator $\hat{\theta}_{n,\Theta_0}$ consistent. If $\Theta_0$ is semi-algebraic and $\theta_0$ a nonsingular interior point of $\Theta_0$, then the likelihood ratio statistic $\lambda_n$ converges to the $\chi_c^2$-distribution with $c = k - \dim_{\theta_0}(\Theta_0)$ degrees of freedom as $n \to \infty$.*

The following is a useful condition for checking the assumption of Theorem 4.3; we use it in Proposition A.4(i) and Theorem 5.1.

LEMMA 4.4. *Let $\Theta_0 = \mathbf{f}(\Gamma)$, where $\Gamma \subseteq \mathbb{R}^d$ is an open semi-algebraic set and $\mathbf{f}$ a polynomial or rational map. If $\theta_0 = \mathbf{f}(\gamma_0) \in \Theta_0$ is nonsingular and the Jacobian of $\mathbf{f}$ at $\gamma_0 \in \Gamma$ of full rank, then $\theta_0$ is in the interior of $\Theta_0$.*

PROOF. The Jacobian being of full rank, there exists a neighborhood $U$ of $\gamma_0$ such that $\mathbf{f}(U)$ is a $d$-dimensional smooth manifold. Since $\theta_0$ is nonsingular, there exists a neighborhood $V$ of $\theta_0$ such that the intersection of $V$ and the Zariski closure $\bar{\Theta}_0$ is also a $d$-dimensional smooth manifold. Since $\mathbf{f}(U) \subseteq \bar{\Theta}_0$, these two manifolds are nested by inclusion. Hence, due to their equal dimension, they must coincide locally. Therefore, in a neighborhood of $\theta_0$, the three sets $\mathbf{f}(U) \subseteq \Theta_0 \subseteq \bar{\Theta}_0$ are equal. It follows that $\theta_0 \in \mathrm{int}(\Theta_0)$. □

Nonsingularity is not necessary for $\chi^2$-asymptotics. For example, suppose $\Theta_0 \in \mathbb{R}^2$ is the union of the two parabolas $y = \pm x^2$ given by the equation $y^2 = x^4$. The origin is a singularity of $\Theta_0$ with tangent cone equal to the $x$-axis. Hence, $\lambda_n \to_d \chi_1^2$ at every point in $\Theta_0$. Removing the part of the parabola passing through the positive orthant, we make the origin a singular boundary point at which $\lambda_n \to_d \chi_1^2$.

4.2. *Algebraic tangent cones and bounds on p-values.* For complicated statistical hypotheses it may be difficult to work out the tangent cone, in which case it is interesting to find sub- and supersets of the tangent cone. The Mahalanobis distances from these sub-/supersets provide distributions that are stochastically larger/smaller than the limiting distribution of the likelihood ratio statistic $\lambda_n$ and thus can be used to bound the asymptotic $p$-value

$$p_\infty(t) = \lim_{n \to \infty} P(\lambda_n > t), \qquad t \geq 0.$$

If a parametrization is available, then the following upper bound is immediate from Lemma 2.4(iv).

LEMMA 4.5. *Suppose $\theta_0 \in \Theta_0 \subseteq \Theta \subseteq \mathbb{R}^k$ is a true parameter point at which the model $\mathcal{P}_\Theta$ is regular and the maximum likelihood estimator $\hat{\theta}_{n,\Theta_0}$ consistent. Let $\Theta_0 = \mathbf{f}(\Gamma)$ be the image of an open semi-algebraic set $\Gamma \subseteq \mathbb{R}^d$ under a polynomial map $\mathbf{f} : \mathbb{R}^d \to \mathbb{R}^k$. Let $J(\gamma)$ be the Jacobian of $\mathbf{f}$ at $\gamma$. If $m$ is the maximum rank of any Jacobian $J(\gamma)$ with $\gamma \in \mathbf{f}^{-1}(\theta_0)$, then $p_\infty(t) \leq P(\chi^2_{k-m} > t)$ for all $t \geq 0$.*

For a lower bound based on a $\chi^2$-distribution one could employ the so-called *Zariski tangent space* given by the kernel of the Jacobian matrix $J(\theta_0)$ that, as in Definition 4.2, is derived from a generating set $\{f_1, \ldots, f_\ell\} \subseteq \mathbb{R}[\mathbf{t}]$ of the ideal $\mathcal{I}(\Theta_0)$. However, if the true parameter $\theta_0$ is a singularity of $\Theta_0$ then the Zariski tangent space is of larger dimension than $\Theta_0$ and does not provide a good local approximation to $\Theta_0$. For instance, the Zariski tangent space at the cusp singularity in Example 1.2 comprises all of $\mathbb{R}^2$ and thus the lower bound is trivially zero because it is computed from a $\chi^2_0$-distribution, which is a point mass at zero.

A better local approximation to $\Theta_0$ is obtained from the algebraic tangent cone defined in (4.3) below. The algebraic tangent cone is sometimes easier to compute than the tangent cone. Gröbner basis methods to automate the computation [7], Section 9.7, are implemented, for example, in `Singular` [15]. We note that in Example 1.2, the algebraic tangent cone at the cusp is equal to the $x$-axis and thus provides a lower $\chi^2_1$-bound for $p_\infty(t)$. A similar phenomenon arises in the feedback model from Example 3.6; see Proposition A.4(iii) in the Appendix.

Let $\theta$ be a point in a semi-algebraic set $\Theta \subseteq \mathbb{R}^k$. For a polynomial $f \in \mathcal{I}(\Theta) \subseteq \mathbb{R}[\mathbf{t}]$ define $f_\theta$ to be the polynomial obtained from $f$ by substituting $t_i + \theta_i$ for each indeterminate $t_i$ appearing in $f$. Write

$$f_\theta = \sum_{h=0}^{j} f_{\theta,h} \tag{4.2}$$

with $f_{\theta,h}$ being a homogeneous polynomial of degree $h$. Define $f_{\theta,\min}$ to be the term $f_{\theta,h}$ that is of smallest degree among all nonzero terms in (4.2). The *algebraic tangent cone* is the real algebraic variety

$$A_\Theta(\theta) = \{\tau \in \mathbb{R}^k \mid f_{\theta,\min}(\tau) = 0 \ \forall f \in \mathcal{I}(\Theta)\}. \tag{4.3}$$

According to the following fact the Mahalanobis distance from the algebraic tangent cone yields a lower bound for $p_\infty(t)$.

LEMMA 4.6. *Let $\theta$ be a point in the semi-algebraic set $\Theta$. Then $T_\Theta(\theta) \subseteq A_\Theta(\theta)$, that is, the algebraic tangent cone contains the tangent cone.*

The algebraic tangent cone $A_\Theta(\theta)$ is a subset of the Zariski tangent space and has dimension equal to the largest dimension of any irreducible component of the Zariski closure $\bar\Theta$ that contains $\theta$. For a polynomial image of an open semi-algebraic set, $\bar\Theta$ is irreducible and the dimension of $A_\Theta(\theta)$ equal to $\dim(\Theta)$. However, the algebraic tangent cone can be of larger dimension than the tangent cone.

The bounds on $p_\infty(t)$ that we discussed in this section are derived from sub- and supersets of the tangent cone. It is noteworthy that if the tangent cone is convex and the limiting distribution a $\chi^2$-mixture, then such bounds can be improved using properties of the mixture weights; compare page 80 in [29]. However, the tangent cones at singularities are generally not convex as can be seen in the example of the feedback model as well as in the factor analysis model that we will study in the remainder of this paper. When studying factor analysis we will employ the following lemma.

LEMMA 4.7. *For $f \in \mathbb{R}[\mathbf{t}]$, define $\Theta$ to be the semi-algebraic set of points $t \in \mathbb{R}^k$ that satisfy $f(t) \geq 0$. If $\theta \in \Theta$ satisfies $f(\theta) = 0$, then the tangent cone $T_\Theta(\theta)$ is contained in the set $\{\tau \in \mathbb{R}^k \mid f_{\theta,\min}(\tau) \geq 0\}$.*

PROOF. Without loss of generality assume that $\theta = 0$ such that $f_\theta = f$. Let $\tau \in T_\Theta(0)$ be the limit of the sequence $(\alpha_n \theta_n)$ with $\alpha_n > 0$ and $\theta_n \in \Theta$ converging to $\theta = 0$. Let $f_{\min} = f_{\theta,\min}$ be of degree $d$. Expanding $f$ as in (4.2) we see that the nonnegative numbers $\alpha_n^d f(\theta_n)$ are equal to $f_{\min}(\alpha_n \theta_n)$ plus a term that converges to zero as $n \to \infty$. Thus, $f_{\min}(\tau) = \lim_{n\to\infty} f_{\min}(\alpha_n \theta_n) \geq 0$. □

**5. Local geometry of the one-factor analysis model.** Assuming zero means to avoid notational overhead, the factor analysis model with $m$ observed variables and $\ell$ hidden factors is the family of multivariate normal distributions $\mathcal{N}_m(0, \Sigma)$ with covariance matrix $\Sigma$ in the set

(5.1) $$F_{m,\ell} = \{\Delta + \Gamma\Gamma^t \mid \Delta \in \mathbb{R}_{\text{pd}}^{m \times m} \text{diagonal}, \Gamma \in \mathbb{R}^{m \times \ell}\}.$$

The set $F_{m,\ell}$ is a semi-algebraic subset of $\mathbb{R}_{\text{sym}}^{m \times m} \simeq \mathbb{R}^{\binom{m+1}{2}}$ with dimension equal to the minimum of $m(\ell + 1) - \binom{\ell}{2}$ and $\binom{m+1}{2}$; see, for example, [10]. This set has singularities and to our knowledge there have been no attempts in the literature to clarify the role these singularities play for asymptotic distribution theory. (Aspects of nonsingular boundary points created by allowing the matrix $\Delta$ in (5.1) to be positive semi-definite have been discussed in [25] and we will not treat these so-called Heywood-cases here.) In this section we derive the tangent cones in factor analysis with $\ell = 1$ factor. The distributional implications are discussed in Section 6.

In the remainder of this section we assume that $m \geq 4$ such that the set $F_{m,1}$ is of dimension $2m < \binom{m+1}{2}$. We begin by describing the ideal $\mathcal{I}(F_{m,1})$. Let $\mathbb{R}[\mathbf{t}]$ be the ring of polynomials in the indeterminates $(t_{ij} \mid 1 \leq i \leq j \leq m)$. Define

$$\mathcal{T}_m = \{t_{ij}t_{gh} - t_{ig}t_{jh}, t_{ih}t_{jg} - t_{ig}t_{jh} \mid 1 \leq i < j < g < h \leq m\} \subset \mathbb{R}[\mathbf{t}].$$

The $2\binom{m}{4}$ quadrics in $\mathcal{T}_m$ are referred to as tetrads in the statistical literature. According to Theorem 16 in [10], the set $\mathcal{T}_m$ generates the ideal $\mathcal{I}(F_{m,1})$.

THEOREM 5.1. *Let $\Sigma = (\sigma_{ij}) \in F_{m,1}$ be the covariance matrix of a distribution in the one-factor analysis model with $m \geq 4$. If there exist at least two nonzero off-diagonal entries $\sigma_{ij}$ and $\sigma_{uv}$ with $i < j$ and $u < v$, then $\Sigma$ is a nonsingular and smooth point of $F_{m,1}$. If at most one off-diagonal entry $\sigma_{ij}$ with $i < j$ is nonzero, then $\Sigma$ is a singular point of $F_{m,1}$. All points $\Sigma \in F_{m,1}$ are of local dimension $\dim_\Sigma(F_{m,1}) = \dim(F_{m,1}) = 2m$.*

PROOF. The claim about local dimension holds because $F_{m,1}$ is the image of a polynomial map; compare Section 3.1. Proposition 32 in [10] states that a matrix $\Sigma = (\sigma_{ij}) \in F_{m,1}$ is a singularity if and only if at most one off-diagonal entry $\sigma_{ij}$ with $i < j$ is nonzero. Let $\mathbf{f}: (0, \infty)^m \times \mathbb{R}^m \to F_{m,1}$ be the parametrization map that sends $(\delta, \Gamma)$ to the matrix $\Delta + \Gamma\Gamma^t \in F_{m,1}$, where $\Delta$ is the diagonal matrix $\mathrm{diag}(\delta)$. In order to show that a nonsingular point $\Sigma = \mathbf{f}(\delta, \Gamma)$ is a smooth point, we check that the $\binom{m+1}{2} \times 2m$-Jacobian of $\mathbf{f}$ at $(\delta, \Gamma)$ is of full rank $2m$; recall Lemma 4.4.

If $\Sigma = \mathbf{f}(\delta, \Gamma)$ has entries $\sigma_{ij} \neq 0$ and $\sigma_{uv} \neq 0$ for two distinct pairs $(i,j)$ and $(u,v)$ with $i < j$ and $u < v$, then $\Gamma = (\gamma_i) \in \mathbb{R}^m$ must have at least three nonzero entries. Without loss of generality, assume that $\gamma_1, \gamma_2, \gamma_3 \neq 0$. Partition the Jacobian matrix of $\mathbf{f}$ by partitioning the columns according to the split between $\delta$ and $\gamma$, and by partitioning the rows into the diagonal and the off-diagonal entries of $\Sigma$. Since $\partial\sigma_{ij}/\partial\delta_k = 0$ if $i < j$, the Jacobian matrix of $\mathbf{f}$ is block-triangular. One of the diagonal blocks, namely the submatrix filled with the partial derivatives $\partial\sigma_{ii}/\partial\delta_j$, is the $m \times m$-identity matrix. Hence, the Jacobian is of full rank $2m$ if and only if the matrix of partial derivatives $\partial\sigma_{ij}/\partial\gamma_k = 0$, $i < j$, is of rank $m$. To see that this rank is indeed $m$, form the $m \times m$-submatrix of partial derivatives $\partial\sigma_{ij}/\partial\gamma_k = 0$ with $(i,j) \in \{(1,2), \ldots, (1,m), (2,3)\}$. This submatrix has determinant equal to $2\gamma_1^{m-2}\gamma_2\gamma_3$ in absolute value. Since $\gamma_1, \gamma_2, \gamma_3 \neq 0$, the determinant is nonzero. □

If $\Sigma \in F_{m,1}$ is a nonsingular point, then the tangent cone $T_{F_{m,1}}(\Sigma)$ is a linear space. At the singularities of $F_{m,1}$ two types of tangent cones arise, which we derive in Lemmas 5.2 and 5.6.

LEMMA 5.2. *If $\Sigma \in F_{m,1}$ is a diagonal matrix, then the tangent cone $T_{F_{m,1}}(\Sigma)$ is equal to the (topological) closure of the set*

$$T_{m,1} = \{\Delta + \Gamma\Gamma^t \mid \Delta \in \mathbb{R}^{m\times m}_{\text{sym}} \text{diagonal}, \Gamma \in \mathbb{R}^m\}.$$

PROOF. In a neighborhood of the origin, the set $F_{m,1} - \Sigma$ is equal to $T_{m,1}$. The set $T_{m,1}$ is a cone such that the claim follows from Lemma 2.4(iii).
□

REMARK 5.3. The cone $T_{m,1}$ is not closed. For example, let

$$\Gamma_n = \left(\frac{1}{n}, n, \frac{1}{n}, 0, \ldots, 0\right)^t \in \mathbb{R}^m,$$

and define $\Delta_n$ to be the diagonal matrix with $i$th diagonal entry equal to the negated square of the $i$th entry of $\Gamma_n$. Then the sequence $(\Delta_n + \Gamma_n\Gamma_n^t)$ converges to a matrix that is zero except for the $(1,2)$- and $(2,3)$-entries that are equal to 1. This limiting matrix is not in $T_{m,1}$.

REMARK 5.4. (i) The result in Lemma 5.2 also carries over to the case of $\ell > 1$ factors. (ii) The algebraic tangent cone $A_{F_{m,1}}(\Sigma)$ at a diagonal matrix $\Sigma$ is simply the real algebraic variety defined by the tetrads $\mathcal{T}_m$.

LEMMA 5.5. *Let $\Sigma \in F_{m,1}$ have exactly one nonzero off-diagonal entry $\sigma_{ij}$ with $i < j$. Then the algebraic tangent cone $A_{F_{m,1}}(\Sigma)$ is equal to the set of matrices $S \in \mathbb{R}^{m\times m}_{\text{sym}}$ that satisfy the following two conditions:* (i) *the $([m] \setminus \{i,j\}) \times ([m] \setminus \{i,j\})$-principal submatrix of $S$ is diagonal, and* (ii) *the rank of the $\{i,j\} \times ([m] \setminus \{i,j\})$-submatrix of $S$ is at most one. Here, $[m] := \{1, \ldots, m\}$.*

PROOF. The set $F_{m,1}$ is the image of the set $(0,\infty)^m \times \mathbb{R}^m$ under a polynomial map. Thus the dimension of $A_{F_{m,1}}(\Sigma)$ is equal to $\dim(F_{m,1}) = 2m$. Without loss of generality we assume that $i = 1$ and $j = 2$.

Let $T = (t_{gh})$ be the symmetric matrix of indeterminates. All $2 \times 2$-minors of the $\{1,2\} \times \{3, \ldots, m\}$-submatrix of $T$ are in the ideal $\mathcal{I}(F_{m,1})$. Since none of these minors involve the indeterminate $t_{12}$, the $\{1,2\} \times \{3, \ldots, m\}$-submatrix of a (symmetric) matrix $S \in A_{F_{m,1}}(\Sigma)$ must have rank at most one.

Let $3 \leq g < h \leq m$. Then the quadric

(5.2) $\quad t_{12}t_{gh} - t_{1g}t_{2h} = \sigma_{12}t_{gh} + (t_{12} - \sigma_{12})t_{gh} - t_{1g}t_{2h}$

is a tetrad in $\mathcal{T}_m$. After substituting $t_{12} + \sigma_{12}$ for $t_{12}$, the lowest degree term on the right-hand side in (5.2) is $\sigma_{12}t_{gh}$. Since $\sigma_{12} \neq 0$, a matrix $S = (s_{ij}) \in A_{F_{m,1}}(\Sigma)$ must have the (off-diagonal) entry $s_{gh} = 0$.

Let $C$ be the set of symmetric matrices for which the $\{1,2\} \times \{3,\ldots,m\}$-submatrix has rank at most one and all off-diagonal entries in the $\{3,\ldots,m\} \times \{3,\ldots,m\}$-submatrix are zero. We have shown that $A_{F_{m,1}}(\Sigma) \subseteq C$. Matrices in $C$ have the $m$ diagonal entries as well as the $(1,2)$-entry unconstrained. Since the set of $2 \times (m-2)$-matrices of rank one has dimension $(m-2)+1 = m-1$, the dimension of $C$ is equal to $2m = \dim(A_{F_{m,1}}(\Sigma))$. The fact that $C$ is an irreducible algebraic variety in $\mathbb{R}^{m \times m}_{\text{sym}}$ now implies $A_{F_{m,1}}(\Sigma) = C$, which is the claim. $\square$

The algebraic tangent cone does not depend on the value of the nonzero off-diagonal entry $\sigma_{ij}$. Unfortunately, this is no longer true for the tangent cone $T_{F_{m,1}}(\Sigma)$.

LEMMA 5.6. *Let $\Sigma \in F_{m,1}$ have exactly one nonzero off-diagonal entry $\sigma_{ij}$ with $i < j$. If $\sigma_{ij} > 0$, then the tangent cone $T_{F_{m,1}}(\Sigma)$ is the set of matrices $S = (s_{gh}) \in A_{F_{m,1}}(\Sigma)$ such that $s_{jg} = \eta \cdot s_{ig}$ for all $g \notin \{i,j\}$ and some $\eta \in [\sigma_{ij}/\sigma_{ii}, \sigma_{jj}/\sigma_{ij}]$. If $\sigma_{ij} < 0$, then the analogue holds with negative multiplier $\eta \in [\sigma_{jj}/\sigma_{ij}, \sigma_{ij}/\sigma_{ii}]$.*

PROOF. Without loss of generality, we assume that $i = 1$ and $j = 2$. Let $\sigma_{12} > 0$ (the case $\sigma_{12} < 0$ is analogous). Denote the set of symmetric matrices claimed to form the tangent cone by $\bar{T}_{F_{m,1}}(\Sigma)$.

We will first show that $T_{F_{m,1}}(\Sigma) \subseteq \bar{T}_{F_{m,1}}(\Sigma)$. We do not change the tangent cone $T_{F_{m,1}}(\Sigma)$ if we restrict $F_{m,1}$ to a neighborhood of $\Sigma$. Hence, we can replace $F_{m,1}$ by $F^\varepsilon_{m,1} = F^\varepsilon_{m,1}(\Sigma)$, which we define to be a neighborhood of $\Sigma$ in $F_{m,1}$ such that $\psi_{12} > 0$ for all $\Psi = (\psi_{ij}) \in F^\varepsilon_{m,1}$. Consider an index $g \geq 3$ and let $\Psi = (\psi_{ij}) = \Delta + \Gamma\Gamma^t \in F^\varepsilon_{m,1}$. If $\Gamma = (\gamma_1,\ldots,\gamma_m)^t$, then $\psi_{12}\psi_{1g}\psi_{2g} = \gamma_1^2\gamma_2^2\gamma_g^2 \geq 0$. It follows that $\psi_{1g}\psi_{2g} \geq 0$ on $F^\varepsilon_{m,1}$. Consequently, $F^\varepsilon_{m,1} = F^{\varepsilon,+}_{m,1} \cup F^{\varepsilon,-}_{m,1}$ with $F^{\varepsilon,+}_{m,1}$ and $F^{\varepsilon,-}_{m,1}$ comprising the matrices $\Psi = (\psi_{ij}) \in F^\varepsilon_{m,1}$ for which $\psi_{1g}, \psi_{2g} \geq 0$ and $\psi_{1g}, \psi_{2g} \leq 0$, respectively. According to Lemma 2.4(i), the tangent cone of $F_{m,1}$ at $\Sigma$ is the union of the two tangent cones of $F^{\varepsilon,+}_{m,1}$ and $F^{\varepsilon,-}_{m,1}$.

Let $\Psi = (\psi_{ij}) = \Delta + \Gamma\Gamma^t \in F^{\varepsilon,+}_{m,1}$ with $\Gamma = (\gamma_1,\ldots,\gamma_m)^t$. Then

$$(5.3) \qquad \psi_{11}\psi_{2g} = (\delta_1 + \gamma_1^2)\gamma_2\gamma_g \geq \gamma_1^2\gamma_2\gamma_g = \psi_{12}\psi_{1g}.$$

Similarly,

$$(5.4) \qquad \psi_{22}\psi_{1g} \geq \psi_{12}\psi_{2g}.$$

Let $S = (s_{ij}) \in T_{F^{\varepsilon,+}_{m,1}}(\Sigma)$. By Lemma 4.7, (5.3) and (5.4) it holds that $s_{1g}, s_{2g} \geq 0$, $\sigma_{11}s_{2g} \geq \sigma_{12}s_{1g}$ and $\sigma_{22}s_{1g} \geq \sigma_{12}s_{2g}$. Thus,

$$(5.5) \quad \text{either} \quad s_{1g} = s_{2g} = 0 \quad \text{or} \quad \left(s_{1g} > 0 \land \frac{s_{2g}}{s_{1g}} \in \left[\frac{\sigma_{12}}{\sigma_{11}}, \frac{\sigma_{22}}{\sigma_{12}}\right]\right),$$

which implies that $s_{2g} = \eta \cdot s_{1g}$ for some $\eta$ as in the claim. Similar consideration of matrices in $F_{m,1}^{\varepsilon,-}$ yields that if $S = (s_{gh}) \in T_{F_{m,1}^{\varepsilon,-}}(\Sigma)$ then $s_{1g}, s_{2g} \leq 0$, $\sigma_{11}s_{2g} \leq \sigma_{12}s_{1g}$ and $\sigma_{22}s_{1g} \leq \sigma_{12}s_{2g}$. This implies an analogue to (5.5), namely,

$$(5.6) \quad \text{either} \quad s_{1g} = s_{2g} = 0 \quad \text{or} \quad \left(s_{1g} < 0 \wedge \frac{s_{2g}}{s_{1g}} \in \left[\frac{\sigma_{12}}{\sigma_{11}}, \frac{\sigma_{22}}{\sigma_{12}}\right]\right),$$

which in turn also implies that $s_{2g} = \eta \cdot s_{1g}$ for some $\eta$ as in the claim. Since the considered index $g \geq 3$ was arbitrary and $T_{F_{m,1}}(\Sigma) \subseteq A_{F_{m,1}}(\Sigma)$, we have proved the inclusion $T_{F_{m,1}}(\Sigma) \subseteq \bar{T}_{F_{m,1}}(\Sigma)$.

In order to show the reverse inclusion, that is, $\bar{T}_{F_{m,1}}(\Sigma) \subseteq T_{F_{m,1}}(\Sigma)$, we write $\Sigma = \Delta_0 + \Gamma_0 \Gamma_0^t$, where $\Delta_0 = \mathrm{diag}(\delta_0)$ is a diagonal and positive definite matrix and $\Gamma_0 = (\gamma_{01}, \gamma_{02}, 0, \ldots, 0)^t \in \mathbb{R}^m$. The pairs $(\gamma_{01}, \gamma_{02})$ that can be used in such a representation of $\Sigma$ satisfy

$$(5.7) \quad \gamma_{01}, \gamma_{02} \neq 0 \quad \text{and} \quad \frac{\gamma_{02}}{\gamma_{01}} \in \left(\frac{\sigma_{12}}{\sigma_{11}}, \frac{\sigma_{22}}{\sigma_{12}}\right),$$

and any value in the interval $(\sigma_{12}/\sigma_{11}, \sigma_{22}/\sigma_{12})$ is possible for their ratio. Let $\mathbf{f} : (0, \infty)^m \times \mathbb{R}^m \to F_{m,1}$ be the parametrization map of $F_{m,1}$; compare the proof of Theorem 5.1. Let $J(\delta, \Gamma)$ be the $\binom{m+1}{2} \times 2m$-Jacobian matrix of $\mathbf{f}$ at $(\delta, \Gamma)$. By Lemma 2.4(iv),

$$\Sigma'(d,c) = J(\delta_0, \Gamma_0)\begin{pmatrix}d \\ c\end{pmatrix} \in \mathbb{R}_{\mathrm{sym}}^{m \times m}, \qquad d, c \in \mathbb{R}^m,$$

is in the tangent cone $T_{F_{m,1}}(\Sigma)$. Let $g, h$ be any two distinct indices in $\{3, \ldots, m\}$. The diagonal entries of $\Sigma'(d,c)$ are

$$\sigma'_{11}(d,c) = d_1 + 2\gamma_{01}c_1, \qquad \sigma'_{22}(d,c) = d_2 + 2\gamma_{02}c_2, \qquad \sigma'_{gg}(d,c) = d_g.$$

Choosing the values $d_i$ appropriately, these diagonal entries may be any real number. The off-diagonal entries of $\Sigma'(d,c)$ are

$$\sigma'_{12}(d,c) = c_1\gamma_{02} + c_2\gamma_{01}, \qquad \sigma'_{1g}(d,c) = \gamma_{01}c_g,$$
$$\sigma'_{gh}(d,c) = 0, \qquad \sigma'_{2g}(d,c) = \gamma_{02}c_g.$$

By appropriate choice of $c_1$ and $c_2$, $\sigma'_{12}(d,c)$ may take on any real value. The entries $\sigma'_{2g}(d,c)$ and $\sigma'_{1g}(d,c)$ are either both zero or both nonzero with their ratio satisfying (5.7). This is equivalent to the existence of a multiplier $\eta$ in the interval in (5.7) such that $\sigma'_{2g}(d,c) = \eta \sigma'_{1g}(d,c)$ for all $g \geq 3$. Therefore, we have shown that any vector in $\bar{T}_{F_{m,1}}(\Sigma)$ for which the multiplier $\eta$ is in the open interval in (5.7) is in the tangent cone $T_{F_{m,1}}(\Sigma)$. However, the tangent cone is a closed set such that the same holds also if $\eta$ is in the closure of the interval in (5.7), which was our claim. $\square$

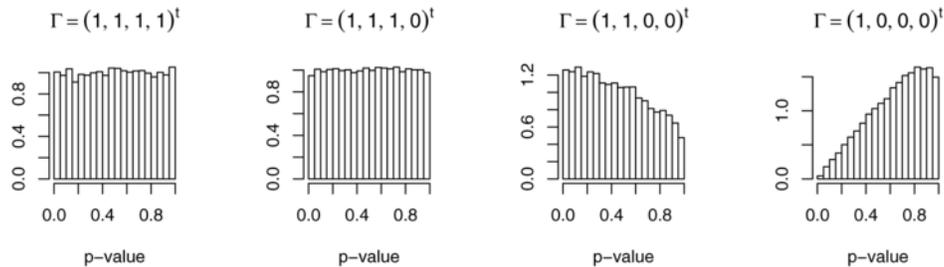

FIG. 3. *Histograms of* 20,000 *simulated p-values for the likelihood ratio test of (6.1) with* $m = 4$ *and sample size* $n = 1000$. *The p-values are computed under* $\chi_2^2$. *The true covariance matrix is equal to* $\Delta + \Gamma\Gamma^t$, *where* $\Delta = \mathrm{diag}(1/3, 1/3, 1/3, 1/3)$ *and* $\Gamma$ *is varied as indicated in the histogram titles. Under these choices pairwise correlations are either zero or equal to* $3/4$.

We remark that the description of the tangent cone in Lemma 5.6 yields a parametrization of the tangent cone. The multiplier $\eta$ in this parametrization is unique unless $s_{jg} = s_{ig} = 0$ for all $g \notin \{i, j\}$.

**6. Likelihood ratio tests in one-factor analysis.** In this section we discuss the limiting distributions for the likelihood ratio statistic $\lambda_n$ in different testing problems involving the factor analysis model with $\ell = 1$ factor.

6.1. *Saturated alternative.* Consider testing the one-factor model against a saturated alternative, that is,

(6.1) $$H_0 : \Sigma \in F_{m,1} \quad \text{vs.} \quad H_1 : \Sigma \notin F_{m,1},$$

where we assume that $m \geq 4$ such that the set $F_{m,1}$ is of positive codimension $\binom{m+1}{2} - 2m$. Statistical software, such as R with command `factanal`, allows one to compute numerically the likelihood ratio statistic $\lambda_n$ for this problem. In such software, $p$-values are computed using the $\chi^2$-distribution with $\binom{m+1}{2} - 2m$ degrees of freedom. Figure 3 shows histograms of simulated $p$-values computed with `factanal`. (Note that `factanal` employs a Bartlett correction.) While the two histograms on the left confirm the expected uniform distribution, this is not the case for the two histograms to the right. It is interesting that the $p$-values for $\Gamma = (1, 1, 0, 0)^t$ tend to be smaller than under a uniform distribution whereas the opposite is true for $\Gamma = (1, 0, 0, 0)^t$.

Figure 3 suggests that there should be at least three different types of limiting distributions for $\lambda_n$. The next result confirms this fact.

THEOREM 6.1. *Let* $\lambda_n$ *be the likelihood ratio statistic for testing (6.1). Assume the true covariance matrix* $\Sigma_0 = (\sigma_{gh})$ *is in* $F_{m,1}$, *and define* $Z \sim \mathcal{N}_{\binom{m+1}{2}}(0, I)$ *to be a standard multivariate normal random vector. When* $n \to \infty$ *it holds that:*

(i) If $\Sigma_0$ has at least two nonzero entries $\sigma_{ij}$ and $\sigma_{uv}$ with $i < j$ and $u < v$, then $\lambda_n$ converges to the $\chi^2$-distribution with $\binom{m+1}{2} - 2m$ degrees of freedom.

(ii) If $\Sigma_0$ has exactly one nonzero off-diagonal entry $\sigma_{ij}$ with $i < j$, then $\lambda_n$ converges to the distribution of the squared Euclidean distance between $Z$ and the (topological) closure of the set of $S = (s_{gh}) \in A_{F_{m,1}}(\Sigma_0)$ such that $s_{jg} = \eta \cdot s_{ig}$ for all $g \notin \{i,j\}$ and some $\eta \geq |\rho_{ij}|/\sqrt{1-\rho_{ij}^2}$. Here, $\rho_{ij} = \sigma_{ij}/\sqrt{\sigma_{ii}\sigma_{jj}}$.

(iii) If $\Sigma_0$ is diagonal, then $\lambda_n$ converges to the distribution of the squared Euclidean distance between $Z$ and the tangent cone $T_{F_{m,1}}(\Sigma_0)$ given in Lemma 5.2.

PROOF. (i) By Theorem 5.1, this is the smooth case and the tangent cone a linear space of the claimed codimension.

(iii) Let $\Sigma_0$ be diagonal. Then the Fisher-information $I(\Sigma_0)$ and its symmetric square root $I(\Sigma_0)^{1/2}$ are diagonal; compare Example 2.2. The diagonal entries of $I(\Sigma_0)$ that are associated with off-diagonal entries $\sigma_{ij}$, $i < j$, factor as
$$I(\Sigma_0)_{ij,ij} = \frac{1}{\sigma_{ii}\sigma_{jj}}.$$
It follows that the tangent cone $T_{F_{m,1}}(\Sigma_0)$ given in Lemma 5.2 is invariant under transformation with $I(\Sigma_0)^{1/2}$. Hence, the claim follows from Theorem 2.6 (recall Lemma 3.3 and Remark 3.4).

(ii) In the remaining case, $\Sigma_0$ has exactly one off-diagonal element, which we assume to be $\sigma_{12} > 0$ (the result for $\sigma_{12} < 0$ is analogous). When listing its rows and columns in the order
$$\underline{11 < 12 < 22} < \underline{13 < 23} < \underline{14 < 24} < \cdots < \underline{1m < 2m} < \underline{33} < \underline{34} < \cdots < \underline{mm}$$
the Fisher-information $I(\Sigma_0)$ is block-diagonal with blocks corresponding to indices that are underlined together. The block for a pair $(1g, 2g)$ with $g \geq 3$ is

(6.2) $$I(\Sigma_0)_{\{1g,2g\}\times\{1g,2g\}} = \frac{1}{\sigma_{gg}}(\Sigma_{12\times 12})^{-1}.$$

Consider the following block-diagonal square root of $I(\Sigma_0)$. For block $\underline{11 < 12 < 22}$ take any square root and for the entries $\underline{33} < \cdots < \underline{mm}$ take the (univariate) square root. For the blocks $\underline{1g < 2g}$ use the Choleski-decomposition of (6.2) to obtain the square root

$$\frac{1}{\sqrt{\sigma_{gg}}}\begin{pmatrix} \frac{1}{\sqrt{\sigma_{11.2}}} & -\frac{\sigma_{12}}{\sigma_{22}\sqrt{\sigma_{11.2}}} \\ 0 & \frac{1}{\sqrt{\sigma_{22}}} \end{pmatrix},$$

where $\sigma_{11.2} = \sigma_{11} - \sigma_{12}^2/\sigma_{22}$. Suppose $\tau = (\tau_{gh})$ is an element of $T_{F_{m,1}}(\Sigma_0)$ for which $\tau_{2g} = \bar{\eta} \cdot \tau_{1g}$ for all $g \geq 3$ and $\bar{\eta} \in (\sigma_{12}/\sigma_{11}, \sigma_{22}/\sigma_{12})$. Under multiplication with the constructed square root of $I(\Sigma_0)$, $\tau$ is mapped to an element $S = (s_{gh})$ of $A_{F_{m,1}}(\Sigma_0)$ for which $s_{2g} = \eta \cdot s_{1g}$ for all $g \geq 3$. The multiplier $\eta$ is equal to

$$(6.3) \qquad \frac{(1/\sqrt{\sigma_{22}})\bar{\eta}}{1/\sqrt{\sigma_{11.2}} - (\sigma_{12}/(\sigma_{22}\sqrt{\sigma_{11.2}}))\bar{\eta}} = \frac{\sqrt{\sigma_{11}\sigma_{22} - \sigma_{12}^2}}{\sigma_{22}/\bar{\eta} - \sigma_{12}}.$$

Therefore,

$$(6.4) \qquad \eta \in \left(\frac{\sigma_{12}}{\sqrt{\sigma_{11}\sigma_{22} - \sigma_{12}^2}}, \infty\right) = \left(\frac{\rho_{12}}{\sqrt{1-\rho_{12}^2}}, \infty\right).$$

We considered $\tau \in T_{F_{m,1}}(\Sigma_0)$ with multiplier $\bar{\eta}$ in the open interval $(\sigma_{12}/\sigma_{11}, \sigma_{22}/\sigma_{12})$. By taking the closure the remaining cases are covered. □

REMARK 6.2. Theorem 12.1 in the seminal paper by Anderson and Rubin [1] gives a sufficient condition for $\chi^2$-asymptotics for the likelihood ratio tests in factor analysis. For the one-factor testing problem (6.1), this theorem states the following. Suppose the true covariance matrix $\Sigma_0 \in F_{m,1}$ is represented as $\Sigma_0 = \Delta + \Gamma\Gamma^t$ with $\Delta$ diagonal and positive definite and $\Gamma \in \mathbb{R}^m$. Then the $\chi^2$-asymptotics from Theorem 6.1(i) hold if the entry-wise (or Hadamard) square of the matrix

$$\Delta - \Gamma(\Gamma^t\Delta^{-1}\Gamma)^{-1}\Gamma^t$$

has nonzero determinant ($\Gamma \neq 0$ is required for this condition to be well defined). We checked that for $m = 4, 5, 6$ this condition is indeed equivalent to requiring two nonzero entries above the diagonal of $\Sigma_0$. However, in the present context, proving Theorem 6.1(i) via Theorem 5.1 seems easier than any attempt to simplify the condition of Anderson and Rubin [1] for the one-factor case.

The distribution described in Theorem 6.1(ii) depends on unknown parameters. This is not the case for the distributional bound obtained from the algebraic tangent cone, for which a nice connection to eigenvalues of Wishart matrices emerges.

THEOREM 6.3. *Let $V$ have a chi-square distribution with $\binom{m-2}{2}$ degrees of freedom and let $W$ be distributed like the smaller of the two eigenvalues of a $2 \times 2$-Wishart matrix with $m - 2$ degrees of freedom and scale parameter the identity matrix $I$. If the true covariance matrix $\Sigma_0 = (\sigma_{gh}) \in F_{m,1}$ has* *exactly one nonzero off-diagonal entry $\sigma_{ij}$ with $i < j$, then the distribution of the squared Mahalanobis distance*

$$\min_{\Sigma \in A_{F_{m,1}}(\Sigma_0)} (Z - \Sigma)^t I(\Sigma_0)(Z - \Sigma), \qquad Z \sim \mathcal{N}_{\binom{m+1}{2}}(0, I(\Sigma_0)^{-1}),$$

*is the distribution of $V + W$, where $V$ and $W$ are independent.*

PROOF. Without loss of generality, assume $\sigma_{12} \neq 0$. We can work with the square root of the Fisher-information $I(\Sigma_0)$ that was used to prove Theorem 6.1(ii). Due to the special block-diagonal structure of $I(\Sigma_0)$, it holds that $A_{F_{m,1}}$ is invariant under transformation with $I(\Sigma_0)^{1/2}$. Thus, the Mahalanobis distance has the same distribution as the Euclidean distance between $\bar{Z} \sim \mathcal{N}_{\binom{m+1}{2}}(0, I)$ and $A_{F_{m,1}}$. The squared Euclidean distance breaks into the sum of

$$V = \sum_{3 \leq g < h \leq m} \bar{Z}_{gh}^2 \sim \chi^2_{\binom{m-2}{2}},$$

and the squared Euclidean distance $W$ between the submatrix $\bar{Z}_{\{1,2\} \times \{3,\ldots,m\}}$ and the set of rank one matrices. The latter distance is equal to the smaller singular value of $\bar{Z}_{\{1,2\} \times \{3,\ldots,m\}}$. The square of this singular value is the smaller eigenvalue $W$ of the $2 \times 2$-Wishart matrix obtained by multiplying $\bar{Z}_{\{1,2\} \times \{3,\ldots,m\}}$ with its transpose. □

Looking back to Theorem 6.1, we see that the $\chi^2$-approximation to the distribution of $\lambda_n$ appears to be valid at almost every covariance matrix in $F_{m,1}$. It is thus tempting to view the singularities as mere theoretical oddities and base inference purely on $\chi^2$-calculations. However, this is problematic because the presence of singularities destroys any possible uniformity of the convergence of $\lambda_n$ to a $\chi^2$-distribution. This can be seen in Figure 4, which shows that the $\chi^2$-approximation becomes more and more inappropriate for smaller and smaller pairwise correlations. A comparison with Figure 3 suggests that this phenomenon is primarily due to the model geometry: small correlations yield points too close to the singular locus of $F_{m,1}$. Indeed the distribution of $\lambda_n$ exhibits features of the limiting distribution from Theorem 6.1(iii); compare the histogram on the far right-hand side in Figure 3.

6.2. *Testing submodels.* In the goodness-of-fit problem (6.1), $\chi^2$-approximations are valid if the true parameter point is far enough away from the singular locus. However, when testing submodels of a factor analysis model, $\chi^2$-approximations may become entirely invalid. We illustrate this for testing the vanishing of some of the components of the parameter vector $\Gamma = (\gamma_1, \ldots, \gamma_m)^t$ defining the covariance matrix $\Sigma = \Delta + \Gamma \Gamma^t \in F_{m,1}$.

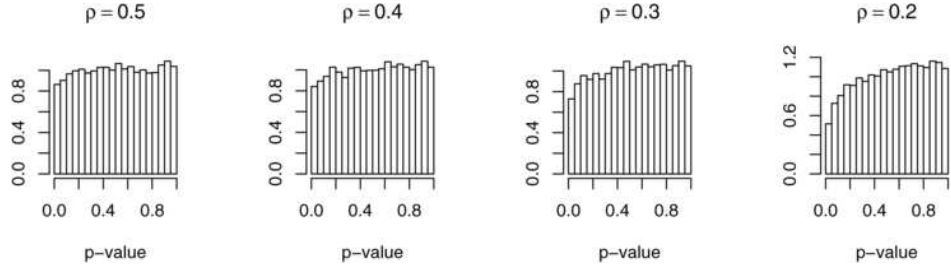

FIG. 4. *Histograms of 20,000 simulated p-values for the likelihood ratio test of (6.1) with $m = 4$ and sample size $n = 50$. The p-values are computed under $\chi_2^2$. The true covariance matrix is a correlation matrix with all off-diagonal entries equal to $\rho$, which is varied as indicated in the histogram titles.*

Let $F_{m,0k}$ be the set of covariance matrices $\Sigma = \Delta + \Gamma\Gamma^t$ in $F_{m,1}$ such that $\Gamma = (\gamma_1, \ldots, \gamma_m)^t$ satisfies that $\gamma_k = \gamma_{k+1} = \cdots = \gamma_m = 0$. Consider testing

(6.5) $\quad H_{0k} : \Sigma \in F_{m,0k} \quad \text{vs.} \quad H_1 : \Sigma \in F_{m,1} \setminus F_{m,0k}.$

Such tests constitute edge exclusion tests in graphical models with one hidden variable; compare, for example, [30]. A positive definite matrix $\Sigma = (\sigma_{ij})$ is in $F_{m,0k}$ if and only if the submatrix $\Sigma_{[k-1]\times[k-1]}$ is in $F_{k-1,1}$ and $\sigma_{ij} = 0$ for all pairs $(i,j) \notin [k-1] \times [k-1]$ with $i \neq j$. Here, $[k-1] = \{1, \ldots, k-1\}$. Hence, the limiting distributions of the likelihood ratio statistic $\lambda_{n,k}$ for testing (6.5) can be determined using Remark 2.8 and Theorems 5.1 and 6.1.

The case $k \geq 4$ is similar to the tests considered in Section 6.1. If $k \geq 4$ and the true covariance matrix $\Sigma_0 \in F_{m,0k}$ cannot be transformed into a matrix in $F_{m,03}$ by permutations of rows and columns, then $\lambda_{n,k}$ converges to a $\chi_{m+k-1}^2$-distribution as $n \to \infty$. At matrices in $F_{m,03}$ (and the possible permutations thereof) nonstandard limiting distributions arise.

The cases $k \leq 3$ are different. If $k = 3$, then there does not exist a true covariance matrix $\Sigma_0 \in F_{m,0k}$ for which $\lambda_{n,k}$ converges to a $\chi^2$-distribution. For $k = 1, 2$, the hypotheses $H_{01}$ and $H_{02}$ are equal because $F_{m,01} = F_{m,02}$ is the set of diagonal covariance matrices. In this case the limiting distribution of $\lambda_{n,k}$ does not depend on $\Sigma_0 \in F_{m,01} = F_{m,02}$. We were not able to connect these distributions to any well-studied distribution but simulations can be used to determine the quantiles of this distribution for a valid (asymptotic) test of $H_{01} = H_{02}$.

When testing $H_{03}$ the limiting distribution of $\lambda_{n,k}$ depends on the correlation $\rho_{12}$. Nevertheless, we have the following corollary to Theorem 6.3; recall Lemma 4.6.

COROLLARY 6.4. *If $\lambda_{n,k}$ is the likelihood ratio statistic for (6.5) with $k = 3$ and the true covariance matrix $\Sigma_0 = (\sigma_{ij})$ is in $F_{m,03}$ with $\sigma_{12} \neq 0$,*

TABLE 1
*Levels of the conservative test that rejects $H_{03}$ if the likelihood ratio statistic exceeds the 95%-quantile of a Wishart eigenvalue distribution. The true covariance matrices are correlation matrices in $F_{m,03}$ with $\rho_{12}$ being varied. Each level was computed in 20,000 simulations*

|  |  | $\rho_{12}$ | | | | | | |
| --- | --- | --- | --- | --- | --- | --- | --- | --- |
|  |  | 0.8 | 0.7 | 0.6 | 0.5 | 0.4 | 0.3 | 0.2 |
| $m=4$ | $n=100$ | 0.027 | 0.028 | 0.034 | 0.034 | 0.036 | 0.042 | 0.048 |
|  | $n=200$ | 0.025 | 0.027 | 0.029 | 0.031 | 0.032 | 0.035 | 0.043 |
|  | $n=500$ | 0.021 | 0.026 | 0.030 | 0.031 | 0.033 | 0.037 | 0.037 |
| $m=8$ | $n=100$ | 0.026 | 0.029 | 0.032 | 0.035 | 0.040 | 0.057 | 0.084 |
|  | $n=200$ | 0.026 | 0.029 | 0.030 | 0.033 | 0.035 | 0.041 | 0.061 |
|  | $n=500$ | 0.023 | 0.027 | 0.027 | 0.031 | 0.032 | 0.035 | 0.040 |

*then*

$$\lim_{n\to\infty} P(\lambda_{n,k} > t) \leq P(W > t),$$

*where $W$ is distributed like the larger of the two eigenvalues of a $2 \times 2$-Wishart matrix with $m-2$ degrees of freedom and scale parameter the identity matrix $I$.*

The algebraic tangent cone calculation that yields Corollary 6.4 thus leads to a simple and conservative test of $H_{03}$: reject the null hypothesis if the observed likelihood ratio statistic is larger than the $1-\alpha$ quantile of the distribution of the eigenvalue $W$. The asymptotic level of this test is provably smaller than $\alpha$ if $\sigma_{12} \neq 0$. Again we point out that there is no uniformity in the convergence to this level and large sample sizes may be required for smaller absolute values of $\sigma_{12}$. Table 1 shows simulated levels for this test using the critical values given in [16]. The increase of the level with $\rho_{12}$ is in agreement with Lemma 5.6 and Theorem 6.1(ii). We note that, if desired, a more powerful yet still asymptotically conservative test can be obtained by relaxing the multiplier $\eta$ in Theorem 6.1(ii) to be in $[0,\infty)$. Critical values for the resulting test of $H_{03}$ could be computed using simulation.

6.3. *Comments on multi-factor analysis.* Factor analysis forms a basic building block for graphical models with hidden variables. As such our study of the one-factor case is of interest for graphical models with one hidden variable. However, in many of its applications factor analysis serves merely as a tool for dimension reduction and the number of factors will typically be much larger than one. While the geometry of models with multiple factors is still largely unknown, the presented theory can offer insights into some of the phenomena encountered by practitioners.

In Section 6.2, we saw that testing the complete independence model, which can be viewed as the case of zero factors, against the one-factor model is a problem for which $\chi^2$-approximations are inappropriate. The generalization of this problem is to test the model with $\ell$ factors against the model with $\ell+1$ factors, that is,

$$(6.6) \qquad H_0 : \Sigma \in F_{m,\ell} \quad \text{vs.} \quad H_1 : \Sigma \in F_{m,\ell+1} \setminus F_{m,\ell}.$$

Simulations such as those in [17] suggest that, regardless of where in the null hypothesis the true distribution is, the likelihood ratio statistic for (6.6) does not follow a $\chi^2$-distribution. Similarly, the limiting distribution when testing

$$(6.7) \qquad H_0 : \Sigma \in F_{m,\ell} \quad \text{vs.} \quad H_1 : \Sigma \notin F_{m,\ell}$$

is not a $\chi^2$-distribution if in fact $\Sigma \in F_{m,k}$ for some $k < \ell$. The algebraic geometrical explanation for these phenomena is that the set $F_{m,\ell+1}$ is singular along its subset $F_{m,\ell}$ ([10], page 491). Note, however, that the set $F_{m,\ell+1}$ has many other more subtle singularities outside $F_{m,\ell}$. These singularities are poorly understood at present.

In the case of $\ell = 1$ factor, singularities arise from independence among observed variables. Consequently, issues with singularities of one-factor models can be avoided if an investigator is free to select variables with pairwise correlations that are large enough for the available sample size. However, problems with singularities are no longer this simple with more than one factor. If $\ell \geq 2$, then detecting whether an estimate of $\Sigma$ is (close to) a singularity of $F_{m,\ell}$ is not a matter of merely gauging whether correlations are different from zero. In the Appendix we illustrate this for the model $F_{5,2}$, which at present is the only model with more than one factor for which the singular locus is known. The details on the algebraic tangent cones of $F_{5,2}$ given in this Appendix show just how complicated the geometry of seemingly simple hidden variable models is.

**7. Conclusion.** We considered likelihood ratio tests of semi-algebraic hypotheses. Using Chernoff's theorem, we showed that under mild probabilistic regularity conditions the large sample limiting distribution of the likelihood ratio statistic always exists. If the true parameter point is a model singularity, then the limiting distribution is determined by the tangent cone. Tangent cones at singularities are generally nonconvex and lead to nonstandard limiting distributions that are different from the mixtures of $\chi^2$-distributions that are often encountered in boundary problems. In fact, singularities can entail arbitrarily complex limiting distributions because any closed semi-algebraic cone of codimension one or larger may occur as tangent cone to a real algebraic variety [12].

Minima of (possibly dependent) $\chi^2$-random variables were seen to be important for locally identifiable models (recall Proposition 3.5). It would be interesting to find good stochastic bounds on the distribution of such minima, which could be used to bound $p$-values when the true parameter point is (close to) a singularity. It seems plausible that bounds could be derived from special constellations of the equi-dimensional tangent spaces that induce the $\chi^2$-variables. Moreover, as pointed out by a referee, the "tube" and "Euler characteristic" methods may be useful for approximating limiting distributions; see [32, 33] and the references therein.

Factor analysis presents interesting examples of models with singularities. Despite its long history and widespread use in practice, these models are far from fully understood. Practical assessment of statistical significance in factor analysis employs $\chi^2$-computations based on the sufficient condition in [1], Theorem 12.1. However, little is known about the structure of the covariance matrices at which $\chi^2$-asymptotics fail and about the nature of the nonstandard limiting distributions. In Sections 5 and 6 we were able to address these problems for factor analysis with one factor.

Factor analysis and all other examples considered in this paper were models for the multivariate normal distribution. Even in this realm there are many other models that could be studied in a similar fashion. For example, more general Gaussian hidden variable models as well as structural equation models could be considered in lieu of factor analysis. But there are also many models for the multinomial distribution that have singularities; see for example [13]. The algebraic geometric techniques presented in this paper provide a unified approach for future work on the impacts of singularities on likelihood ratio tests in different classes of nonsmooth models.

A key feature of the $\chi^2$-theory for smooth models is that the limiting distribution is pivotal, that is, does not depend on where in the null hypothesis the true parameter point is located. In some nonstandard problems, such as testing submodels of the one-factor model (Section 6.2) this pivotality is preserved (or at least stochastic bounds are pivotal). However, our computations for the simplest nontrivial two-factor model suggest that even for testing problems involving only slightly more general and still seemingly simple hidden variable models, the limiting distribution will depend on unknown nuisance parameters. One possible approach to circumventing this problem is to design bootstrap procedures. While this is a topic beyond the scope of this paper, we expect the algebraic framework layed out here will be helpful for investigating asymptotic correctness of bootstrap tests in the presence of singularities.

## APPENDIX A: DETAILS ON THE FEEDBACK MODEL

Here we provide details on the feedback model $\mathcal{P}_{\Theta_0}$ from Example 3.6; the same notation is used.

**A.1. Parameter identifiability.** Instead of studying the parameterization map **f** directly, we work with precision matrices $\Sigma^{-1}$. For $\beta \in \mathbb{R}^5$ and $\kappa \in (0, \infty)^4$, let **g** be the map from $(\beta, \kappa)$ to the inverse of the covariance matrix $\mathbf{f}(\beta, \omega)$ with $\omega_i = 1/\kappa_i$. The map $\mathbf{g} : \Gamma \subsetneq \mathbb{R}^5 \times (0,\infty)^4 \to \mathbb{R}^{4 \times 4}_{\text{sym}}$ is polynomial with $\mathbf{g}(\beta, \kappa)$ equal to

$$\begin{pmatrix} \kappa_1 + \beta_{21}^2 \kappa_2 + \beta_{31}^2 \kappa_3 & \beta_{32}\beta_{31}\kappa_3 - \beta_{21}\kappa_2 & -\beta_{31}\kappa_3 & \beta_{24}\beta_{21}\kappa_2 \\ & \kappa_2 + \beta_{32}^2 \kappa_3 & -\beta_{32}\kappa_3 & -\beta_{24}\kappa_2 \\ & & \kappa_3 + \beta_{43}^2 \kappa_4 & -\beta_{43}\kappa_4 \\ & & & \kappa_4 + \beta_{24}^2 \kappa_2 \end{pmatrix}.$$

(A.1)

Computations in `Maple` and `Singular` [15] yield the following results.

The Jacobian of **g** has an $8 \times 8$-minor equal to a product of powers of the $\kappa_i$. Hence, its rank is 8 or 9 for all $(\beta, \kappa) \in \Gamma$. Computing the radical of the ideal of $9 \times 9$-minors, we see that the rank is 9 unless (3.2) holds. In order to investigate identifiability of $\mathcal{P}_{\Theta_0}$ we perform computations for solving the polynomial equations $\mathbf{g}(\bar{\beta}, \bar{\kappa}) = \mathbf{g}(\beta, \kappa)$ for $(\bar{\beta}, \bar{\kappa})$. We find the following structure:

LEMMA A.1. *The model $\mathcal{P}_{\Theta_0}$ is globally identifiable at $(\beta, \kappa)$ if and only if:*

(i) $\beta_{31} + \beta_{32}\beta_{21} \neq 0$, *or*

(ii) $\beta_{31} + \beta_{32}\beta_{21} = 0$ *and at least one of the parameters $\beta_{32}$, $\beta_{43}$, $\beta_{24}$ is zero, or*

(iii) $\beta_{31} + \beta_{32}\beta_{21} = 0$ *and* $\beta_{32}\beta_{43}\beta_{24} = -1$.

If $\mathcal{P}_{\Theta_0}$ is not globally identifiable at $(\beta, \kappa)$ and $\Sigma^{-1} = \mathbf{g}(\beta, \kappa)$, then $\mathbf{g}^{-1}(\Sigma^{-1}) = \{(\beta, \kappa), (\bar{\beta}, \bar{\kappa})\}$ is of cardinality two; compare (3.3). The nontrivial element $(\bar{\beta}, \bar{\kappa})$ is a rational function of $(\beta, \kappa)$. It holds that $\bar{\beta}_{21} = \beta_{21}$ and $\bar{\kappa}_1 = \kappa_1$. For $i \in \{2, 3, 4\}$,

$$(A.2) \quad \bar{\beta}_{i\,i-1} = \frac{\kappa_{i-1}\kappa_{i+1} + \kappa_i \kappa_{i+1} \beta_{i\,i-1}^2 + \kappa_i \kappa_{i-1} \beta_{i\,i-1}^2 \beta_{i-1\,i+1}^2}{\beta_{i\,i-1}(\kappa_i \kappa_{i+1} + \kappa_i \kappa_{i-1}\beta_{i-1\,i+1}^2 + \kappa_{i-1}\kappa_{i+1}\beta_{i-1\,i+1}^2 \beta_{i+1\,i}^2)}$$

and

$$(A.3) \qquad\qquad \bar{\kappa}_i = \kappa_i \beta_{i\,i-1} / \bar{\beta}_{i\,i-1}.$$

The expressions in (A.2) and (A.3) can be read literally for $i = 3$; for $i = 2, 4$ the indices $i \pm 1$ are to be read modulo 3 such that $4 + 1 \equiv 2$ and $2 - 1 \equiv 4$. Finally, the remaining component $\bar{\beta}_{31}$ is equal to $-\bar{\beta}_{32}\bar{\beta}_{21}$.

LEMMA A.2. *If $\mathcal{P}_{\Theta_0}$ is locally identifiable at $(\beta, \kappa) \in \Gamma$, then the map **g** is proper at the precision matrix $\Sigma^{-1} = (\sigma^{ij}) = \mathbf{g}(\beta, \kappa)$.*

PROOF. Let $V \subseteq \mathbb{R}^{4 \times 4}_{\text{sym}}$ be a compact neighborhood of $\Sigma^{-1}$ in $\mathbf{g}(\Gamma)$. We assume in particular that $V$ is bounded away from the boundary of the cone of positive definite matrices such that the closure of $\mathbf{g}^{-1}(V)$ is contained in $\Gamma$ (recall that points in $\Gamma$ satisfy $\beta_{43}\beta_{32}\beta_{24} \neq 1$). By the local identifiability assumption, none of the values of $\beta_{32}$, $\beta_{43}$ and $\beta_{24}$ are zero, which implies that the off-diagonal entries $\sigma^{ij}$ are nonzero if $i, j \geq 2$. We may thus assume that for all precision matrices $S = (s_{ij}) \in V$, the diagonal entries $s_{ii}$ as well as the absolute values of off-diagonal entries $s_{ij}$ with $i, j \geq 2$ are in some finite interval $[m, M]$ with $0 < m < M < \infty$.

Suppose $S = (s_{ij}) = \mathbf{g}(\beta, \kappa) \in V$. Since $\kappa_i \leq s_{ii}$ it follows that all the components of $\kappa$ are in the interval $[0, M]$. For the components of $\beta$ we have the inequalities

$$\beta_{31}^2 m^2 \leq \beta_{31}^2 s_{23}^2 = \beta_{31}^2 \beta_{32}^2 \kappa_3^2 \leq s_{11} s_{22} \leq M^2,$$
$$\beta_{21}^2 m^2 \leq \beta_{21}^2 s_{24}^2 = \beta_{21}^2 \beta_{24}^2 \kappa_2^2 \leq s_{11} s_{44} \leq M^2$$

and

$$|\beta_{ij}| m \leq |\beta_{ij} s_{ij}| = \beta_{ij}^2 \kappa_i \leq s_{ii} \leq M, \qquad (i, j) \in \{(4, 3), (3, 2), (2, 4)\}.$$

The absolute values of the components of $\beta$ are thus all contained in the interval $[0, M/m]$. Hence, $(\beta, \kappa)$ is in the compactum $[-M/m, M/m]^5 \times [0, M]^4$. $\square$

Proposition 3.5 does not apply to case (ii) in Lemma A.1 because if one of $\beta_{32}$, $\beta_{43}$, $\beta_{24}$ is zero, then $\mathbf{g}$ is not proper at $\mathbf{g}(\beta, \kappa)$. This can be seen in equation (A.2) where $|\bar{\beta}_{i\,i-1}| \to \infty$ if $\beta_{i\,i-1} \to 0$. In case (iii) of Lemma A.1, Proposition 3.5 does not apply because the rank of the Jacobian of $\mathbf{g}$ drops from 9 to 8.

**A.2. Applying algebraic techniques.** Example 3.6 was concerned with local identifiability. In order to get an understanding of the globally identifiable cases the techniques from Section 4 are useful.

LEMMA A.3. *A covariance matrix $\Sigma = \mathbf{f}(\beta, \omega) \in \Theta_0$ is a singularity if and only if $\beta_{31} + \beta_{32}\beta_{21} = 0$, in which case the algebraic tangent cone $A_{\Theta_0}(\Sigma)$ is the union of the two hyperplanes with the normal vectors $\eta$ and $\bar{\eta}$ from (3.4) and (3.5).*

PROOF. Employing the technique of implicitization ([7], Chapter 3) and the software `Singular` [15], we compute the Zariski closure $\bar{\Theta}_0$ that is found to be given by the vanishing of the irreducible polynomial

$$f = \sigma_{13}\sigma_{14}^3\sigma_{23}^2 - 2\sigma_{13}^2\sigma_{14}^2\sigma_{23}\sigma_{24}$$

$$+ \sigma_{13}^3\sigma_{14}\sigma_{24}^2 - \sigma_{12}\sigma_{14}^3\sigma_{23}\sigma_{33} + \sigma_{12}\sigma_{13}\sigma_{14}^2\sigma_{24}\sigma_{33}$$
$$+ \sigma_{11}\sigma_{14}^2\sigma_{23}\sigma_{24}\sigma_{33} - \sigma_{11}\sigma_{13}\sigma_{14}\sigma_{24}^2\sigma_{33}$$
$$+ \sigma_{12}^2\sigma_{14}^2\sigma_{33}\sigma_{34} - \sigma_{11}\sigma_{14}^2\sigma_{22}\sigma_{33}\sigma_{34}$$
$$- \sigma_{12}^2\sigma_{13}\sigma_{14}\sigma_{34}^2 + \sigma_{11}\sigma_{13}\sigma_{14}\sigma_{22}\sigma_{34}^2$$
$$+ \sigma_{12}\sigma_{13}^2\sigma_{14}\sigma_{23}\sigma_{44} - \sigma_{11}\sigma_{13}\sigma_{14}\sigma_{23}^2\sigma_{44}$$
$$- \sigma_{12}\sigma_{13}^3\sigma_{24}\sigma_{44} + \sigma_{11}\sigma_{13}^2\sigma_{23}\sigma_{24}\sigma_{44}$$
$$- \sigma_{12}^2\sigma_{13}\sigma_{14}\sigma_{33}\sigma_{44} + \sigma_{11}\sigma_{13}\sigma_{14}\sigma_{22}\sigma_{33}\sigma_{44}$$
$$+ \sigma_{12}^2\sigma_{13}^2\sigma_{34}\sigma_{44} - \sigma_{11}\sigma_{13}^2\sigma_{22}\sigma_{34}\sigma_{44}.$$

We can also compute the singularities of $\Theta_0$, which are the matrices $\Sigma = (\sigma_{ij})$ with $\sigma_{13} = \sigma_{14} = 0$, or equivalently the matrices $\mathbf{f}(\beta, \omega)$ with $\beta_{31} + \beta_{32}\beta_{21} = 0$.

Since the ideal $\mathcal{I}(\Theta_0)$ is generated by $f$, the algebraic tangent cone at a singularity $\Sigma = \mathbf{f}(\beta, \omega)$ is determined by the polynomial $f_{\Sigma,\min} \in \mathbb{R}[t]$, which factorizes as

$$(\beta_{43} \cdot t_{13} - t_{14})[(\beta_{32}^2\beta_{43}^2\omega_2 + \beta_{43}^2\omega_3 + \omega_4) \cdot t_{13}$$
$$- (\omega_3 + \beta_{32}^2\omega_2 + \beta_{32}^2\beta_{24}^2\omega_4)\beta_{43} \cdot t_{14}].$$

The linear forms in the factorization correspond to the vectors in (3.4) and (3.5). □

The next proposition summarizes what we know about the globally identifiable cases from Lemma A.1. Note that case (ii) is an example where the parametrization is globally identifiable with full rank Jacobian, but where a nonstandard limiting distribution arises for the likelihood ratio statistic.

PROPOSITION A.4. *Let $\lambda_n$ be the likelihood ratio statistic for testing (1.1). Let $\Sigma_0 = \mathbf{f}(\beta, \omega)$ be the true covariance matrix. Suppose $n \to \infty$.*

(i) *If $\beta_{31} + \beta_{32}\beta_{21} \neq 0$, then $\lambda_n \to \chi_1^2$.*

(ii) *If $\beta_{31} + \beta_{32}\beta_{21} = 0$ and at least one of the parameters $\beta_{32}, \beta_{43}, \beta_{24}$ is zero, then $\lambda_n$ converges to the distribution of a minimum of two $\chi_1^2$-random variables, which as in Example 3.6 is determined by the cosine $\rho$ in (3.6).*

(iii) *If $\beta_{31} + \beta_{32}\beta_{21} = 0$ and $\beta_{32}\beta_{43}\beta_{24} = -1$, then the asymptotic p-values for the likelihood ratio test can be bounded as $P(\chi_1^2 > t) \leq p_\infty(t) \leq P(\chi_2^2 > t)$.*

PROOF. (i) This is the smooth case, which follows from (3.2), Lemmas A.3 and 4.4, and Theorem 4.3.

(ii) Under the assumed conditions on $\beta$, the algebraic tangent cone is the union of two distinct hyperplanes because $\eta \neq \bar{\eta}$. The hyperplane given by $\eta$ corresponds to the span of the columns of the Jacobian of $\mathbf{f}$ at $(\beta, \omega)$.

The other hyperplane comprises tangent vectors obtained from diverging sequences $(\bar{\beta}, \bar{\omega})$ such that $\mathbf{f}(\bar{\beta}, \bar{\omega}) \to \Sigma_0$. This can be done by plugging the expressions in (A.2) and (A.3) into the Jacobian of $\mathbf{f}$ (note that $\omega_i = 1/\kappa_i$) and computing the column span. Hence, $T_{\Theta_0}(\Sigma_0) = A_{\Theta_0}(\Sigma_0)$ is of the same form as for the locally identifiable case discussed in Example 3.6. (If $\beta_{43} = \beta_{24} = 0$ or $\beta_{43} = \beta_{32} = 0$, then $\rho = 0$ implies independence of the two $\chi_1^2$-random variables of which the minimum is taken.)

(iii) In this case, the normals $\eta$ and $\bar{\eta}$ are proportional to each other and the two associated hyperplanes cone coincide. Hence, $A_{\Theta_0}(\Sigma_0)$ is a hyperplane and $P(\chi_1^2 > t) \leq p_\infty(t)$. Since the rank of the Jacobian of $\mathbf{f}$ at $(\beta, \omega)$ is equal to 8, the upper bound on $p_\infty(t)$ follows from Lemma 4.5. $\square$

The tangent cone $T_{\Theta_0}(\Sigma_0)$ in case (iii) of Proposition A.4 seems difficult to obtain. However, based on examination of the degree 3-terms in the equation that defines $\Theta_0 - \Sigma_0$ at the origin, we believe that, similarly to Example 1.2, $T_{\Theta_0}(\Sigma_0) \subsetneq A_{\Theta_0}(\Sigma_0)$.

## APPENDIX B: THE SIMPLEST TWO-FACTOR MODEL

In this appendix we discuss the geometry of $F_{5,2}$, the covariance matrix parameter space of factor analysis with 5 observed variables and 2 factors. The Zariski closure of $F_{5,2}$ is the hypersurface defined by the vanishing of the *pentad*

$$t_{12}t_{13}t_{24}t_{35}t_{45} - t_{12}t_{13}t_{25}t_{34}t_{45} - t_{12}t_{14}t_{23}t_{35}t_{45} + t_{12}t_{14}t_{25}t_{34}t_{35}$$
$$+ t_{12}t_{15}t_{23}t_{34}t_{45} - t_{12}t_{15}t_{24}t_{34}t_{35} + t_{13}t_{14}t_{23}t_{25}t_{45} - t_{13}t_{14}t_{24}t_{25}t_{35}$$
$$- t_{13}t_{15}t_{23}t_{24}t_{45} + t_{13}t_{15}t_{24}t_{25}t_{34} - t_{14}t_{15}t_{23}t_{25}t_{34} + t_{14}t_{15}t_{23}t_{24}t_{35}.$$

Finding the tangent cones of $F_{5,2}$ is an open problem but we can compute the algebraic tangent cones of this hypersurface.

The singularities of $F_{5,2}$ are of two types ([10], Example 33). First, there are the symmetric matrices with a row (and column) that is off-diagonally zero. Consider the matrices

$$(\text{B.1}) \qquad \Sigma = \begin{pmatrix} \sigma_{11} & 0 & 0 & 0 & 0 \\ 0 & \sigma_{22} & \sigma_{23} & \sigma_{24} & \sigma_{25} \\ 0 & \sigma_{23} & \sigma_{33} & \sigma_{34} & \sigma_{35} \\ 0 & \sigma_{24} & \sigma_{34} & \sigma_{44} & \sigma_{45} \\ 0 & \sigma_{25} & \sigma_{35} & \sigma_{45} & \sigma_{55} \end{pmatrix}$$

with first row and column off-diagonally zero as a representative set. For almost all singularities $\Sigma = (\sigma_{ij})$ of form (B.1), the algebraic tangent cone $A_{F_{5,2}}(\Sigma)$ is the irreducible real algebraic variety given by the quadratic polynomial

$$\sigma_{45}(\sigma_{24}\sigma_{35} - \sigma_{34}\sigma_{25})t_{12}t_{13} - \sigma_{35}(\sigma_{23}\sigma_{45} - \sigma_{34}\sigma_{25})t_{12}t_{14} \pm \cdots$$

obtained by replacing the indeterminates $t_{gh}$ with $2 \leq g < h \leq 5$ in the pentad by $\sigma_{gh}$. However, at special matrices of form (B.1) the algebraic tangent cone may be of degree 3 or larger. This occurs if the submatrix $\Sigma_{\{2,\ldots,5\} \times \{2,\ldots,5\}}$ satisfies all its tetrads or has an off-diagonal $2 \times 2$-submatrix that is zero. Degree 4 occurs if precisely one entry above the diagonal of $\Sigma$ is nonzero. If $\Sigma$ is diagonal, then $A_{F_{5,2}}(\Sigma)$ is the pentad hypersurface itself.

The second type of singularities of $F_{5,2}$ is given by symmetric matrices that satisfy all those tetrads that do not involve some given off-diagonal entry. As a representative set, consider the matrices that satisfy all those tetrads that do not involve $\sigma_{12}$. We note that these matrices can be parametrized as

$$(\text{B.2}) \quad \Sigma = \begin{pmatrix} \delta_{11} & 0 & 0 & 0 & 0 \\ 0 & \delta_{22} & 0 & 0 & 0 \\ 0 & 0 & \delta_{33} & 0 & 0 \\ 0 & 0 & 0 & \delta_{44} & 0 \\ 0 & 0 & 0 & 0 & \delta_{55} \end{pmatrix} + \begin{pmatrix} \gamma_{11} & \gamma_{12} \\ \gamma_{21} & \gamma_{22} \\ \gamma_{31} & 0 \\ \gamma_{41} & 0 \\ \gamma_{51} & 0 \end{pmatrix} \begin{pmatrix} \gamma_{11} & \gamma_{12} \\ \gamma_{21} & \gamma_{22} \\ \gamma_{31} & 0 \\ \gamma_{41} & 0 \\ \gamma_{51} & 0 \end{pmatrix}^t.$$

For almost all singularities $\Sigma = (\sigma_{ij})$ of form (B.2), the algebraic tangent cone $A_{F_{5,2}}(\Sigma)$ is the irreducible real algebraic variety defined by a quadratic polynomial with 18 terms

$$\sigma_{15}\sigma_{45}t_{23}t_{34} + \sigma_{34}\sigma_{45}t_{23}t_{15} + \sigma_{23}\sigma_{34}t_{45}t_{15} - \sigma_{25}\sigma_{45}t_{34}t_{13} - \sigma_{15}\sigma_{35}t_{34}t_{24} \pm \cdots$$

that are obtained from the six pentad monomials of the form $t_{12}t_{2i_3}t_{i_3i_4}t_{i_4i_5}t_{1i_5}$ by dropping $t_{12}$, and replacing $t_{2i_3}t_{i_3i_4}$, $t_{i_3i_4}t_{i_4i_5}$ or $t_{i_4i_5}t_{1i_5}$ by the corresponding expression in $\sigma_{gh}$. The degree of $A_{F_{5,2}}(\Sigma)$ is larger than two if (i) $\gamma_{12} = 0$ or $\gamma_{22} = 0$ such that $\Sigma \in F_{5,1}$, or (ii) two or more of the coefficients $\gamma_{31}$, $\gamma_{41}$ and $\gamma_{51}$ are zero which leads to at least two off-diagonally zero rows (and columns) in $\Sigma$, or (iii) $\gamma_{11}$, $\gamma_{21}$ and at least one coefficient among $\gamma_{31}$, $\gamma_{41}$ and $\gamma_{51}$ are zero. In case (iii) the matrix becomes block-diagonal; for example, if $\gamma_{11} = \gamma_{21} = \gamma_{31} = 0$ then $\Sigma = \text{diag}(\Sigma_{12 \times 12}, \sigma_{33}, \Sigma_{45 \times 45})$. As for the singularities of the first type, the algebraic tangent cone admits degree 4 if $\Sigma$ has precisely one nonzero off-diagonal entry and degree 5 if $\Sigma$ is diagonal.

For a generic one-factor matrix $\Sigma \in F_{5,1}$, the cone $A_{F_{5,2}}(\Sigma)$ is given by a cubic polynomial with 60 terms

$$\sigma_{15}\sigma_{45}t_{12}t_{23}t_{34} + \sigma_{12}\sigma_{15}t_{23}t_{34}t_{45} + \sigma_{34}\sigma_{45}t_{12}t_{23}t_{15}$$
$$+ \sigma_{23}\sigma_{34}t_{12}t_{45}t_{15} + \sigma_{12}\sigma_{23}t_{34}t_{45}t_{15} - \sigma_{25}\sigma_{45}t_{12}t_{34}t_{13}$$
$$- \sigma_{12}\sigma_{25}t_{34}t_{45}t_{13} - \sigma_{24}\sigma_{45}t_{23}t_{15}t_{13} \pm \cdots.$$

The terms are obtained by choosing one of the twelve monomials $t_{i_1i_2}t_{i_2i_3}t_{i_3i_4}t_{i_4i_5}t_{i_5i_1}$ in the pentad and one of the five distinct indices $i_j$, and replacing $t_{i_ji_{j+1}}t_{i_{j+1}i_{j+2}}$ by $\sigma_{i_ji_{j+1}}\sigma_{i_{j+1}i_{j+2}}$. Here the additions $j+1$ and $j+2$ are modulo 5.

The above algebraic tangent cones depend on the numerical values of the entries of a singularity. Since the tangent cones themselves might even be more diverse, as was the case with one-factor analysis, we can expect the likelihood ratio statistic for testing the two-factor model to admit many different limiting distributions.

**Acknowledgments.** The author would like to thank two referees and an Associate Editor for their helpful comments on an earlier draft of this paper.

Department of Statistics
University of Chicago
Chicago, Illinois 60637
USA
E-mail: drton@galton.uchicago.edu